\documentclass[preprint,3p]{elsarticle}
\usepackage{amssymb,amsmath,amsthm,indentfirst}
\usepackage{multicol}
\usepackage{graphicx}
\DeclareMathOperator{\Rea}{Re}
\DeclareMathOperator{\Ima}{Im}
\renewcommand{\leq}{\leqslant}
\renewcommand{\geq}{\geqslant}
\newproof{pf}{Proof}

\newtheorem{proposition}{Proposition}
\newtheorem{lemma}{Lemma}

\journal{Applied Numerical Mathematics}

\begin{document}

\begin{frontmatter}

\title{A Numerov-Crank-Nicolson-Strang scheme with discrete transparent boundary conditions for the Schr\"odinger equation on a semi-infinite strip}

\author{Alexander Zlotnik \corref{cor2}\fnref{fn1,fn2}}
\ead{azlotnik2008@gmail.com}
\author{Alla Romanova \corref{cor1}\fnref{fn1}}
\ead{avromm1@gmail.com}

\cortext[cor2]{Principal corresponding author}
\cortext[cor1]{Corresponding author}
\fntext[fn1]{Department of Higher Mathematics at Faculty of Economics,
National Research University Higher School of Economics,
Myasnitskaya 20, 101000 Moscow, Russia.}
\fntext[fn2]{Department of Mathematical Modelling,
National Research University Moscow Power Engineering Institute,
Krasnokazarmennaya 14, 111250 Moscow, Russia.}
\begin{abstract}
\noindent We consider an initial-boundary value problem for a 2D time-dependent Schr\"odinger equation on a semi-infinite strip.
For the Numerov-Crank-Nicolson finite-difference scheme with discrete transparent boundary conditions, the Strang-type splitting with respect to the potential is applied.
For the resulting method, the uniqueness of a solution and the uniform in time $L^2$-stability (in particular, $L^2$-conservativeness) are proved.
Due to the splitting, an effective direct algorithm using FFT in the direction perpendicular to the strip is developed to implement the splitting method for general potential.
Numerical results on the tunnel effect for smooth and rectangular barriers together with the practical error analysis on refining meshes are included as well.
\end{abstract}


\end{frontmatter}
\par MSC[2010] 65M06, 65M12, 35Q41.
\section{\large Introduction}
\label{sect1}
The time-dependent Schr\"odinger equation with several variables is important in quantum mechanics, atomic and nuclear physics, wave physics, microelectronics and nanotechnologies, etc. Often it should be solved in unbounded space domains.
\par Several approaches were developed and studied to deal with problems of such kind, in particular, see \cite{AABES08,ABM04,AES03,DiM97,Sch02a,SZ04}.
One of them exploits the so-called discrete transparent boundary conditions (TBCs) at artificial boundaries \cite{AES03,EA01}.
Its advantages are the complete absence of spurious reflections, reliable computational stability, clear mathematical background and corresponding rigorous stability theory.
\par Concerning finite-difference schemes in several space variables with the discrete TBCs, the standard Crank-Nicolson scheme in the case of an infinite (or semi-infinite) strip was studied in detail in \cite{AES03,DZ06,DZ07}; the higher order Numerov-Crank-Nicolson scheme was considered in \cite{SA08}, and
a general family of schemes was treated in \cite{ZZ11,IZ11}.
But all the schemes are implicit, and solving of specific complex systems of linear algebraic equations is required at each time level. Efficient methods to solve similar systems are well developed by the moment in real situation but not the complex one.
\par The splitting technique is widely used to simplify numerical solving of the time-dependent Schr\"o\-dinger and related equations, in particular, see \cite{BM00,CMR13}, \cite{GX11}-\cite{NT09}, \cite{TY10}.
The Strang-type splitting with respect to the potential (but not to the space directions) has been recently applied to the Crank-Nicolson scheme with the discrete TBC in \cite{DZZ13}.
\par It is well-known that higher order methods are often able to reduce computational costs significantly. In this paper, we apply the Strang-type splitting in potential to the Numerov-Crank-Nicolson scheme and get a higher order in space method with the discrete TBC.
To study its stability and construct the discrete TBC, we first consider the splitting Numerov-Crank-Nicolson scheme on an infinite space mesh in the semi-infinite strip. Its uniform in time $L^2$-stability together with the mass conservation law are proved.
\par The discrete TBC allows to restrict in a rigorous sense the solution of the latter scheme to a finite space mesh. Our form of the discrete TBC is different from (though  equivalent to) one for the original Numerov-Crank-Nicolson scheme in \cite{SA08}. Its presented derivation based on results in \cite{DZZ09} is more direct and shorter (we skip an intermediate step dealing with the computationally unstable form of the discrete TBC). Notice that the explicit form of the discrete TBC is non-local and involves the discrete convolution in time as well as the discrete Fourier expansion in direction $y$ perpendicular to the strip.
We prove both the uniqueness of solution exploiting some results from \cite{DZZ09,ZL10} and the uniform in time $L^2$-stability of the resulting method. In particular, it is $L^2$-conservative.
Due to the splitting, an effective direct algorithm using FFT in $y$ is developed to implement the method (for general potential).
\par The corresponding numerical results on the tunnel effect for two barriers (smooth and rectangular) are presented.
They are accompanied by the practical error analysis using the numerical solutions on refining meshes in all space and time directions.
In the case of the rectangular barrier, the results improve those from \cite{DZZ13} by taking coarser space mesh and thus reducing computational costs due to higher order in space of the method. Notice that the average of the discontinuous rectangular barrier is essential.
\section{\large The Schr\"odinger equation in a semi-infinite strip
 and its approximations of higher order in space}
\label{sect4}
\setcounter{equation}{0}
\setcounter{proposition}{0}
\setcounter{theorem}{0}
\setcounter{lemma}{0}
\setcounter{corollary}{0}
\setcounter{remark}{0}
We consider the 2D time-dependent Schr\"odinger equation
\begin{equation}
 i\hbar\frac{\partial\psi}{\partial t}
 =-\frac{\hbar^{\,2}}{2m_0}\,\Delta \psi+V\psi\ \
 \text{for}\ \ (x,y)\in\Pi_{\infty},\ \ t>0,
\label{a1}
\end{equation}
where $\Delta$ is the 2D Laplace operator,
$\Pi_{\infty}:=(0,\infty)\times (0,Y)$ is a semi-infinite strip and $V(x,y)$ is a given real potential.
Also $i$ is the imaginary unit, $\hbar>0$ and $m_0>0$ are physical constants, and the unknown wave function $\psi=\psi(x,y,t)$ is complex-valued.
Below we use an abbreviation $c_{\hbar}:=\frac{\hbar^{\,2}}{2m_0}$.
\par We impose the following boundary condition, condition at infinity and initial condition
\begin{gather}
 \psi(\cdot,t)|_{\partial\Pi_{\infty}}=0,\ \
 \|\psi(x,\cdot,t)\|_{L^2(0,Y)}\to 0\ \
 \text{as}\ \ x\to\infty,
 \ \ \text{for any}\ \ t>0,
\label{a2}\\[1mm]
 \psi|_{t=0}=\psi^0(x,y)\ \ \text{in}\ \ \Pi_{\infty}.
 \label{a3}
\end{gather}
We also assume that $V(x,y)$ is constant and $\psi_0(x,y)$ vanishes when $x$ is sufficiently large:
\begin{equation}
 V(x,y)=V_{\infty},\ \
 \psi^0(x,y)=0\ \ \text{for}\ \ (x,y)\in [X_0,\infty)\times (0,Y),
\label{a4}
\end{equation}
 for some $X_0>0$.
\par We introduce a uniform mesh $\overline{\omega}_{hx}$ in $x$ with nodes $x_j=jh_x$, $j\geq 0$, and the step $h_x=\frac{X}{J}$,
where $X>X_0$ and $h_x\leq X-X_0$.
We introduce a uniform mesh $\overline{\omega}_{hy}$ in $y$ with nodes $y_k=kh_y$, $0\leq k\leq K$, and the step $h_y=\frac{Y}{K}$.
Let $\omega_{hy}:=\{y_k\}_{k=1}^{K-1}$.
\par We define the product mesh $\overline{\omega}_{h,\infty}:=\overline{\omega}_{hx}\times\overline{\omega}_{hy}$, its internal part $\omega_{h,\infty}:=\left\{(x_j,y_k);\,j\geq 1\right.$, $\left.1\leq k\leq K-1\right\}$ and the boundary $\Gamma_{h,\infty}:=\overline{\omega}_{h,\infty}\backslash\omega_{h,\infty}$. Hereafter $h=(h_x,h_y)$.
\par We define the backward, forward and central difference quotients as well as the average in $x$
\begin{gather}
 \bar{\partial}_xW_j:= \frac{W_j-W_{j-1}}{h_x},\ \
 \partial_xW_j:= \frac{W_{j+1}-W_{j}}{h_x},\ \
 \overset{\circ}{\partial}_x W_j:= \frac{W_{j+1}-W_{j-1}}{2h_x},
\nonumber\\[1mm]
s_{\theta x}W_j:= \theta W_{j-1}+(1-2\theta)W_{j}+\theta W_{j+1}
=(I+\theta h_x^2\partial_x\bar{\partial}_x)W_j
\label{b1}
\end{gather}
with the parameter $\theta$, where $I$ is the unit operator (it will be often omitted below).
The average becomes the Numerov one for $\theta=\frac{1}{12}$
\begin{equation}
s_{Nx}W_j:=\frac{1}{12}\, W_{j-1}+\frac{5}{6}\, W_{j}+\frac{1}{12}\, W_{j+1}
=\Bigl(I+\frac{h_x^2}{12}\, \partial_x\bar{\partial}_x\Bigr)W_j.
\label{b2}
\end{equation}
We also define the corresponding difference quotients and the Numerov average in $y$
\begin{gather*}
 \bar{\partial}_yU_k:= \frac{U_k-U_{k-1}}{h_y},\ \
 \partial_yU_k:= \frac{U_{k+1}-U_{k}}{h_y},
\\
s_{Ny}U_k:= \frac{1}{12}\, U_{k-1}+\frac{5}{6}\,U_{k}+\frac{1}{12}\,U_{k+1}
=\Bigl(I+\frac{h_y^2}{12}\,\partial_y\bar{\partial}_y\Bigr)U_k.
\end{gather*}
\par We introduce also a non-uniform mesh $\overline{\omega}^{\,\tau}$ in time on $[0,\infty)$ with nodes
$0=t_0< t_1<\dots< t_m<\dots$, where $t_m\rightarrow\infty$ as $m\rightarrow\infty$, and steps $\tau_m=t_m-t_{m-1}$.
Let $\omega^{\tau}:=\overline{\omega}^{\,\tau}\backslash \{0\}$ and $\tau_{\max}=\sup_{m\geq 1} \tau_m$.
We define the backward difference quotient, the symmetric average and the backward shift in time
\[
 \bar{\partial}_t Z^m= \frac{Z^m-Z^{m-1}}{\tau_m},\ \
 \overline{s}_t Z^m= \frac{Z^{m-1}+Z^m}{2},\ \
 {\check Z}^m=Z^{m-1}.
\]
\par The Numerov-type discretizations of the Laplace operator and the 2D average are given by
\begin{gather*}
\Delta_{hN}:=s_{Ny}\partial_x\bar{\partial}_x+s_{Nx}\partial_y\bar{\partial}_y
=\Delta_h+\frac{h_x^2+h_y^2}{12}\,\partial_x\bar{\partial}_x\partial_y\bar{\partial}_y,
\ \
s_N:=I+\frac{h_x^2}{12}\,\partial_x\bar{\partial}_x+\frac{h_y^2}{12}\,\partial_y\bar{\partial}_y,
\end{gather*}
where $\Delta_h:=\partial_x\bar{\partial}_x+\partial_y\bar{\partial}_y$ is the simplest approximation of the 2D Laplace operator.
\par We begin with a known discretization for the Schr\"odinger equation (\ref{a1}) which is of the Numerov type in space and two-level symmetric in time (i.e. of the Crank-Nicolson type) given by
\begin{equation}
i\hbar s_N\bar{\partial}_t \Psi
=-c_{\hbar}\Delta_{hN}\overline{s}_t \Psi
+s_N\left( V\overline{s}_t \Psi\right)
\ \ \mbox{on}\ \ \omega_{h,\infty}\times\omega^{\tau}.
\label{c2}
\end{equation}
The corresponding approximation error
\[
i\hbar s_N\bar{\partial}_t \psi
+c_{\hbar}\Delta_{hN}\overline{s}_t \psi
-s_N\left( V\overline{s}_t \psi\right)
=O\left(\tau_{\max}^2+|h|^4\right)
\]
is the 2nd order in $\tau_{\max}$ and of higher 4th order in $|h|$, for $\psi$ smooth enough.
\par We supplement the discrete equation (\ref{c2}) with the boundary and initial conditions
\begin{equation}
\left. \Psi\right|_{\Gamma_{h,\infty}\times\omega^{\tau}}=0,
\ \ \Psi^0=\Psi^0_h\ \ \mbox{on}\ \ \overline{\omega}_{h,\infty}.
\label{f2}
\end{equation}
Hereafter the compatibility condition $\left. \Psi^0_h\right|_{\Gamma_{h,\infty}}=0$ is assumed to be valid.
\par
We further apply the Strang-type splitting in the potential to the scheme (\ref{c2}), (\ref{f2}) and get the following three-step scheme
\begin{gather}
 i\hbar\,\frac{\breve{\Psi}^m-\Psi^{m-1}}{\tau_m/2}
 =\Delta V\frac{\breve{\Psi}^m+\Psi^{m-1}}{2}\ \
 \text{on}\ \ \omega_{h,\infty},
\label{e1}\\[1mm]
 i\hbar s_N\frac{\widetilde{\Psi}^m-\breve{\Psi}^m}{\tau_m}
 =-c_{\hbar}{\Delta}_{hN}\frac{\widetilde{\Psi}^m+\breve{\Psi}^m}{2}
 +s_N\Bigl(\tilde{V}\frac{\widetilde{\Psi}^m+\breve{\Psi}^m}{2}\Bigr)
\ \ \text{on}\ \ \omega_{h,\infty},
\label{e2}\\[1mm]
 i\hbar\,\frac{\Psi^m-\widetilde{\Psi}^m}{\tau_m/2}
 =\Delta V\frac{\Psi^m+\widetilde{\Psi}^m}{2}\ \
\ \ \text{on}\ \ \omega_{h,\infty},
\label{e3}
\end{gather}
with the boundary and initial conditions
\begin{gather}
 \breve{\Psi}^m|_{\Gamma_h,\infty}=0,\ \ \widetilde{\Psi}^m|_{\Gamma_h,\infty}=0,\ \
 \Psi^m|_{\Gamma_h,\infty}=0,
\label{e4}\\[1mm]
 \Psi^0=\Psi^0_h\ \ \text{on}\ \ \overline{\omega}_{h,\infty},
\label{e5}
\end{gather}
for any $m\geq 1$, where $\Delta V:=V-\tilde{V}$ and the auxiliary 1D potential $\tilde{V}=\tilde{V}(x)$ satisfies $\tilde{V}(x)=V_\infty$ for $x\geq X_0$. In the simplest case, $\tilde{V}(x)=V_\infty$ (but a non-constant $\tilde{V}$ is required when extending the results to the case of an infinite strip with different values $V_{\pm\infty}$ of $V(x,y)$ when $x\to\pm\infty$). \par The construction of this splitting is similar to the case of the standard scheme without averages \cite{DZZ13}. Note that we have omitted operators $s_N$ arising in the course of the splitting on both sides of (\ref{e1}) and (\ref{e3}). Clearly equations \eqref{e1} and \eqref{e3}
are reduced to the explicit formulas
\begin{equation}
 \breve{\Psi}^m=\mathcal{E}^m\Psi^{m-1},\ \
 \Psi^m=\mathcal{E}^m\widetilde{\Psi}^m,\ \ \text{with}\ \
 \mathcal{E}^m
 :=\frac{1-i\displaystyle{\frac{\tau_m}{4\hbar}}\, \Delta V}
        {1+i\displaystyle{\frac{\tau_m}{4\hbar}}\, \Delta V}.
\label{e6}
\end{equation}
The main equation (\ref{e2}) is similar to the original one (\ref{c2}) on the time level $m$ but it is simplified by
replacing the coefficient $V(x,y)$ by $\tilde{V}(x)$. The functions $\breve{\Psi}$ and $\widetilde{\Psi}$ are auxiliary unknowns while $\Psi$ is the main one.
\par From (\ref{e6}) we get immediately
\begin{equation}
 |\breve{\Psi}^m|=|\Psi^{m-1}|,\ \ |\Psi^m|=|\widetilde{\Psi}^m|\ \ \mbox{on}\ \ \overline{\omega}_{h,\infty};
\label{ff1}
\end{equation}
moreover, since $\Delta V_{j,k}=0$ for $j\geq J-1$, we simply have
\begin{equation}
 \breve{\Psi}^m_{j,k}=\Psi^{m-1}_{j,k},\ \ \Psi^m_{j,k}=\widetilde{\Psi}^m_{j,k}\ \ \mbox{for}\
 \ j\geq J-1.
\label{ff2}
\end{equation}
\par This splitting modifies the scheme (\ref{c2}), (\ref{f2}) only in time and is symmetric in time due to steps (\ref{e1}) and (\ref{e3}). Thus, concerning the approximation error, it diminishes neither the 4th order in $|h|$ nor the 2nd order in $\tau_{\max}$. This can be checked more formally as follows. Eliminating the auxiliary functions $\breve{\Psi}$ and $\widetilde{\Psi}$ from (\ref{e2}) with the help of (\ref{e6}) leads to the following equation for $\Psi$
\begin{equation}
 i\hbar s_N\, \frac{{\mathcal E}^{-1}\Psi-{\mathcal E} \check{\Psi}}{\tau}
=-c_{\hbar} {\Delta}_{hN}\frac{{\mathcal E}^{-1}\Psi+{\mathcal E}\check{\Psi}}{2}+s_N\Bigl(\tilde{V}\frac{{\mathcal E}^{-1}\Psi+{\mathcal E}\check{\Psi}}{2}\Bigr).
\label{f3}
\end{equation}
Note that ${\mathcal E}^{-1}={\mathcal E}^*$ and
\begin{gather*}
\frac{{\mathcal E}^*\Psi-{\mathcal E} \check{\Psi}}{\tau}
=(\Rea {\mathcal E})\bar{\partial}_t \Psi
-i\, \frac{2}{\tau}\, (\Ima {\mathcal E})\overline{s}_t \Psi,
\\
\frac{{\mathcal E}^*\Psi+{\mathcal E} \check{\Psi}}{2}
=(\Rea {\mathcal E})\overline{s}_t \Psi
-i\, \frac{\tau}{2}\, (\Ima {\mathcal E})\bar{\partial}_t \Psi.
\end{gather*}
Consequently we can rewrite equation (\ref{f3}) as follows
\begin{gather*}
i\hbar s_N
 \left[ (\Rea {\mathcal E})\bar{\partial}_t \Psi \right]
=-c_{\hbar} {\Delta}_{hN}\left[(\Rea {\mathcal E})\overline{s}_t \Psi
-i\, \frac{\tau}{2}\, (\Ima {\mathcal E})\bar{\partial}_t \Psi \right]
\\[1mm]
+s_N\Bigl[\Bigl((\Rea {\mathcal E})\tilde{V}-\frac{2\hbar}{\tau}\, \Ima {\mathcal E}\Bigr)\overline{s}_t \Psi
-i\, \frac{\tau}{2}\, (\Ima {\mathcal E})\tilde{V}\bar{\partial}_t \Psi
 \Bigr].
\end{gather*}
After some calculations, this formulation implies that the approximation error of the discrete equation (\ref{f3}) differs from the original one (\ref{c2}) by a term of order $O\left(\tau_{\max}^2\right)$ and thus they are of the same order $O\left(\tau_{\max}^2+|h|^4\right)$ (in particular, note that
$\Rea {\mathcal E}\to 1$ and $-\frac{2\hbar}{\tau}\, \Ima {\mathcal E}\to \Delta V$ as $\tau_{\max}\to 0$).
On the other hand, now checking the same approximation order imposes additional smoothness restrictions on $V$.
\section{\large Stability of the splitting higher order scheme on an infinite space mesh}
\setcounter{equation}{0}
\setcounter{proposition}{0}
\setcounter{theorem}{0}
\setcounter{lemma}{0}
\setcounter{corollary}{0}
\setcounter{remark}{0}
Let $H_h$ be a Hilbert space of mesh functions $W$: $\overline{\omega}_{h,\infty}\to {\mathbb C}$ such that
\[
\left. W\right|_{\Gamma_{h,\infty}}=0\ \ \text{and}\ \
\sum_{j=1}^{\infty} \sum_{k=1}^{K-1}\left|W_{j,k}\right|^2<\infty
\]
endowed with the following mesh counterpart of the inner product on $L^2(\Pi_{\infty})$
\[
\left(U,W\right)_{H_h}
:=\sum_{j=1}^{\infty} \sum_{k=1}^{K-1}U_{j,k}W_{j,k}^* h_xh_y.
\]
\par In order to study the stability in more detail, we treat a non-homogeneous version of (\ref{e2})
\begin{equation}
i\hbar s_N\, \frac{\widetilde{\Psi}^m-\breve{\Psi}^m}{\tau_m}
 =-c_{\hbar}{\Delta}_{hN}\frac{\widetilde{\Psi}^m+\breve{\Psi}^m}{2}
 +s_N\Bigl(\tilde{V}\frac{\widetilde{\Psi}^m+\breve{\Psi}^m}{2} \Bigr)
+F^m
\ \ \text{on}\ \ \omega_{h,\infty}.
\label{h1}
\end{equation}
We use only the first assumption (\ref{a4}) in the section.

\begin{proposition}
\label{ph}
Let $\Psi_h^0, F^m\in H_h$ for any $m\geq 1$.
Then there exists a unique solution to the splitting scheme \eqref{e1}, \eqref{h1}, \eqref{e3}-\eqref{e5}
such that $\Psi^m\in H_h$ for any $m\geq 0$, and the following $L^2$-stability bound holds
\begin{equation}
 \max_{0\leq m\leq M}\|\Psi^m\|_{H_h}
 \leq\|\Psi^0_h\|_{H_h}
 +\frac{6}{\hbar}\sum_{m=1}^M\left\|F^m\right\|_{H_h}\tau_m\ \ \text{for any}\ \ M\geq 1.
\label{h2}
\end{equation}
\par Moreover, in the particular case $F=0$, the following mass conservation law holds
\begin{equation}
 \|\Psi^m\|_{H_h}^2=
 \|\Psi^0_h\|_{H_h}^2
 \ \ \text{for any}\ \ m\geq 1.
\label{h3}
\end{equation}
\end{proposition}
\begin{pf}
We first rewrite the main equation (\ref{h1}) as a suitable operator equation in $H_h$. We set $\Lambda_xW:=-\partial_x\bar{\partial}_xW$ and $\Lambda_yW:=-\partial_y\bar{\partial}_yW$ on $\omega_{h,\infty}$
and $\Lambda_xW:=\Lambda_yW:=0$ on $\Gamma_{h,\infty}$. Then the operators
\begin{gather*}
 \Lambda_x,\,\ \Lambda_y,\ \
 s_{Nx}=I-\frac{h_x^2}{12}\,\Lambda_x,\ \
 s_{Ny}=I-\frac{h_y^2}{12}\,\Lambda_y,\,\
 s_N=I-\frac{h_x^2}{12}\,\Lambda_x-\frac{h_y^2}{12}\,\Lambda_y
\end{gather*}
and $-{\Delta}_{hN}=s_{Ny}\Lambda_x+s_{Nx}\Lambda_y$
are bounded and self-adjoint in $H_h$. Moreover, the following inequalities hold
\begin{gather}
\left(s_{Nx}W,W\right)_{H_h}\geq \frac{2}{3} \left\|W\right\|_{H_h}^2,\ \
\left(s_{Ny}W,W\right)_{H_h}\geq \frac{2}{3} \left\|W\right\|_{H_h}^2\ \
\mbox{for any}\ W\in H_h
\label{h3a}
\end{gather}
(see \cite{DZZ09} in the particular case of a space average with parameter $\theta=\frac{1}{12}$) and
\[
\left(s_NW,W\right)_{H_h}\geq \frac{1}{3} \left\|W\right\|_{H_h}^2\ \
\mbox{for any}\ W\in H_h
\]
since $s_N=s_{Nx}+s_{Ny}-I$.
Therefore the inverse operator $s_N^{-1}$ exists and is bounded:
\begin{equation}
 \left\|s_{Nx}^{-1}\right\|_{{\mathcal L}(H_h)}\leq 3.
\label{i1}
\end{equation}
Thus we can consider (\ref{h1}) as an operator equation in $H_h$. Following \cite{SA08}, we apply
$s_N^{-1}$
to it and obtain
\begin{equation}
i\hbar\,\frac{\widetilde{\Psi}^m-\breve{\Psi}^m}{\tau_m}
 =A_h\,\frac{\widetilde{\Psi}^m+\breve{\Psi}^m}{2}+s_N^{-1}F^m
\ \ \text{in}\ \ H_h,
\label{i2}
\end{equation}
where
\[
A_h:=c_{\hbar}\left[\Lambda_x+\Lambda_y
+s_N^{-1}\left(\frac{h_x^2}{12}\,\Lambda_x^2+\frac{h_y^2}{12}\,\Lambda_y^2\right)\right]+\tilde{V}.
\]
Since $s_N$ commutes with $\Lambda_x$ and $\Lambda_y$, so $s_N^{-1}$ does with $\Lambda_x^2$ and $\Lambda_y^2$.
Consequently, $A_h$ is a bounded self-adjoint operator in $H_h$.
\par We rewrite equation (\ref{i2}) in another form
\[
\left(I+i\frac{\tau_m}{2\hbar}\, A_h\right)\widetilde{\Psi}^m
=B^m
:=
\left(I-i\frac{\tau_m}{2\hbar}\, A_h\right)\breve{\Psi}^m
-i\frac{\tau_m}{\hbar}\,s_N^{-1}F^m.
\]
Since the operator $I+i\frac{\tau_m}{2\hbar}\, A_h$ is invertible, the equation has a unique solution $\widetilde{\Psi}^m\in H_h$ provided that
$\breve{\Psi}^m,F^m\in H_h$. This implies the existence of a unique solution of the splitting scheme such that $\Psi^m\in H_h$ for any $m\geq 0$.
\par We can now apply a technique from \cite{DZZ13}. First note that the pointwise equalities (\ref{ff1}) imply
\begin{equation}
 \|\breve{\Psi}^m\|_{H_h}=\|\Psi^{m-1}\|_{H_h},\ \ \|\Psi^m\|_{H_h}=\|\widetilde{\Psi}^m\|_{H_h}.
\label{j1}
\end{equation}
Multiplying the operator equation (\ref{i2}) by $\frac{\widetilde{\Psi}^m+\breve{\Psi}^m}{2}$, separating the imaginary part of the result and using the property $A_h=A_h^*$, we get
\[
\frac{\hbar}{2\tau_m} \left(
\|\widetilde{\Psi}^m\|_{H_h}^2-\|\breve{\Psi}^m\|_{H_h}^2
\right)
=
\Ima
\Bigl(
 s_N^{-1}F^m,
\frac{\widetilde{\Psi}^m+\breve{\Psi}^m}{2}
\Bigr)_{H_h}.
\]
Applying equalities (\ref{j1}), multiplying both sides by $\frac{2\tau_m}{\hbar}$ and summing up the result over $m=1,\dots,M$, we obtain
\begin{equation}
\|\Psi^M\|_{H_h}^2
=\|\Psi^0_h\|_{H_h}^2
+\frac{2}{\hbar}
\sum_{m=1}^M
\Ima
\Bigl(
  s_N^{-1}F^m,
\frac{\widetilde{\Psi}^m+\breve{\Psi}^m}{2}
\Bigr)_{H_h} \tau_m.
\label{j2}
\end{equation}
Due to (\ref{i1}) and (\ref{j1}) we further have
\begin{gather*}
 \|\Psi^m\|_{H_h}^2
 \leq\|\Psi^0\|_{H_h}^2
 +\frac{2}{\hbar}\sum_{m=1}^M
 \|s_N^{-1}F^m\|_{H_h}
 \frac{\|\widetilde{\Psi}^m\|_{H_h}+\|\breve{\Psi}^m\|_{H_h}}{2}\, \tau_m
\\[1mm]
 \leq
 \|\Psi^0\|_{H_h}^2
 +\frac{6}{\hbar}
 \sum_{m=1}^M \|F^m\|_{H_h} \max_{0\leq m\leq M}\|\Psi^m\|_{H_h} \tau_m.
\end{gather*}
This bound directly implies (\ref{h2}).
The conservation law (\ref{h3}) follows from (\ref{j2}).
\end{pf}
\section{\large The splitting higher order scheme on a finite space mesh}
\setcounter{equation}{0}
\setcounter{proposition}{0}
\setcounter{theorem}{0}
\setcounter{lemma}{0}
\setcounter{corollary}{0}
\setcounter{remark}{0}
The splitting scheme (\ref{e1})-(\ref{e5}) is not practically implementable because of the infinite number of unknowns on each time level.
We intend to restrict its solution to a finite space mesh $\overline{\omega}_h:=\{(x_j,y_k)\in \overline{\omega}_{h,\infty};\, 0\leq j\leq J\}$.
Let $\omega_h:=\{(x_j,y_k)\in \omega_{h,\infty};\, 1\leq j\leq J-1$, $1\leq k\leq K-1\}$ be its internal part and
$\partial\omega_h=\overline{\omega}_h\backslash \omega_h$ be its boundary.
We also set $\Gamma_{1h}:=\{(x_J,y_k);\, 1\leq k\leq K-1\}$
and $\Gamma_h=\partial\omega_h\backslash \Gamma_{1h}$.
\par By definition, \textit{the discrete transparent boundary condition} (TBC) is a boundary condition on $\Gamma_{1h}$ which allows one to accomplish the restriction.
\par We need the well-known direct and inverse discrete Fourier sine transforms in $y$
\begin{gather*}
Z^{(q)}=\left({\mathcal F}_yZ\right)^{(q)}:= \frac{2}{K} \sum_{k=1}^{K-1} Z_j \sin \frac{\pi q k}{K},\ \ 1\leq q\leq K-1,
\\[1mm]
Z_k=\left({\mathcal F}_y^{-1}Z^{(\cdot)}\right)_k:= \sum_{q=1}^{K-1} Z^{(q)} \sin \frac{\pi q k}{K},\ \ 1\leq k\leq K-1.
\end{gather*}
The corresponding eigenvalues of the operators $-\partial_y\bar{\partial}_y$ and $s_{Ny}$ are
\[
 \lambda_q=
 \left(\frac{2}{h_y}\sin \frac{\pi q h_y}{2Y}\right)^2\ \ \mbox{and}\ \
 \sigma_q=1-\frac{h_y^2}{12}\lambda_q
 =1-\frac{1}{3}\sin^2 \frac{\pi q h_y}{2Y}
\in \Bigl(\frac23,1\Bigr),\ \ 1\leq q\leq K-1.
\]
\par Given a function $W$: $\overline{\omega}_h\to {\mathbb C}$, denote by $W_j$ its restriction to $\{x_j\}\times\overline{\omega}_{hy}$; in particular, $W_J$ is its trace on $\Gamma_{1h}$.
Let ${\bf \Psi}^m_J$ be the vector function ${\bf \Psi}^m_J=\{\Psi_J^0,\dots,\Psi_J^m\}$.
Let also
\[
(R*Q)^m:=\sum_{p=0}^m R^pQ^{m-p},\ \ m\geq 0,
\]
be the discrete convolution product of the mesh functions $R,Q$: $\overline{\omega}^{\,\tau}\to {\mathbb C}$.
\par
Assume that the time mesh be uniform with a step $\tau>0$ below.
\begin{proposition}
\label{pk}
Let $\Psi_h^0=0$ and $ F^m=0$ on
$\omega_{h,\infty}\backslash\omega_h$ for any $m\geq 1$ and
$\bigl.\Psi_{h}^0\bigr|_{j=J-1}=0$.

The solution to the splitting scheme \eqref{e1}, \eqref{h1}, \eqref{e3}-\eqref{e5}
such that $\Psi^m\in H_h$ for any $m\geq 0$ satisfies the following three-step splitting scheme on the finite space mesh $\overline{\omega}_h$
\begin{gather}
 i\hbar\, \frac{\breve{\Psi}^m-\Psi^{m-1}}{\tau/2}
 =\Delta V\frac{\breve{\Psi}^m+\Psi^{m-1}}{2}\ \
 \text{on}\ \ \omega_h\cup \Gamma_{1h},
\label{l1}\\[1mm]
 i\hbar s_N\,\frac{\widetilde{\Psi}^m-\breve{\Psi}^m}{\tau}
 =-c_{\hbar}{\Delta}_{hN}\frac{\widetilde{\Psi}^m+\breve{\Psi}^m}{2}
 +s_N\Bigl(\tilde{V}\frac{\widetilde{\Psi}^m+\breve{\Psi}^m}{2}\Bigr)
 +F^m
\ \ \text{on}\ \ \omega_h,
\label{l2}\\[1mm]
 i\hbar\,\frac{\Psi^m-\widetilde{\Psi}^m}{\tau/2}
 =\Delta V\frac{\Psi^m+\widetilde{\Psi}^m}{2}\ \
 \text{on}\ \ \omega_h\cup \Gamma_{1h},
\label{l3}
\end{gather}
with the boundary and initial conditions
\begin{gather}
 \breve{\Psi}^m|_{\Gamma_h}=0,\ \ \widetilde{\Psi}^m|_{\Gamma_h}=0,\ \
 \Psi^m|_{\Gamma_h}=0,
\label{l4}\\[1mm]
 \Psi^0=\Psi^0_h\ \ \text{on}\ \ \overline{\omega}_h  ,
\label{l6}
\end{gather}
as well as the discrete TBC for the Fourier coefficients of the solution
\begin{gather}
\Bigl.\Bigl[
c_{\hbar,q}\bar{\partial}_x\,\frac{\widetilde{\Psi}^{m(q)}+\breve{\Psi}^{m(q)}}{2}
-h_xs_{\theta_q x}^-
\Bigl(
i\hbar\,\frac{\widetilde{\Psi}^{m(q)}-\breve{\Psi}^{m(q)}}{\tau}
-V_{\infty,q}\,\frac{\breve{\Psi}^{m(q)}+\widetilde{\Psi}^{m(q)}}{2}
\Bigr)
\Bigr]\Bigr|_{j=J}
\nonumber\\[1mm]
=c_{\hbar,q}
\bigl(R_q *{\widetilde\Psi}^{(q)}_J\bigr)^m
\label{l5m}
\end{gather}
for any $m\geq 1$ and $1\leq q\leq K-1$, where
\begin{equation}
 c_{\hbar,q}:=c_{\hbar}\Bigl[1+\Bigl(\frac{h_xh_y\lambda_q}{12\sigma_q}\Bigr)^2\Bigr],\ \
 \theta_q:=\frac{1}{12\sigma_q}\in \Bigl(\frac{1}{12},\frac18\Bigr),\ \
 V_{\infty,q}:= V_{\infty}+c_{\hbar}\frac{\lambda_q}{\sigma_q}
\label{l51}
\end{equation}
and $s_{\theta x}^{-}W_j:=\theta W_{j-1}+(\frac12-\theta) W_j$. Hereafter $Z^{m(q)}=(Z^m)^{(q)}$.
\par The discrete convolution kernel in (\ref{l5m})
can be computed by the recurrent formulas
\begin{gather}
 R_q^0=c_{1q},\ \ R_q^1=-c_{1q}\varkappa_q \mu_q,
\label{rr1}\\[1mm]
 R_q^m=\frac{2m-3}{m}\, \varkappa_q \mu_q R^{m-1}_q-\frac{m-3}{m}\, \varkappa_q^2 R_q^{m-2}\ \ \mbox{for}\ \ m\geq 2,
\label{rr2}
\end{gather}
with the coefficients defined by
\begin{gather*}
 c_{1q}=
 -\frac{|\alpha_q|^{1/2}}{2}\,e^{-i(\arg \alpha_q)/2},\ \
 \varkappa_q=-e^{i\arg \alpha_q},\ \
 \arg \alpha_q\in(0,2\pi),\ \
 \mu_q=\frac{\beta_q}{|\alpha_q|}\in(-1,1),
\\[1mm]
 \alpha_q=2a_q+(1-4\theta_q)h_x^2a_q^2\neq 0,\ \
 \beta_q=2\Rea a_q+(1-4\theta_q)h_x^2|a_q|^2,\ \
 a_q=\frac{V_{\infty,q}}{2c_{\hbar,q}}+i\, \frac{\hbar}{\tau c_{\hbar,q}}.
\end{gather*}
\end{proposition}
\begin{pf}
Clearly it suffices only to derive the discrete TBC (\ref{l5m}) for the solution of the splitting scheme (\ref{e1}), (\ref{h1}), (\ref{e3})-(\ref{e5}) under the above assumptions on $\Psi^0_h$ and $F$.
\par Due to property (\ref{f2}), equations (\ref{e1}), (\ref{h1}) and (\ref{e3}) are reduced
on $(\omega_{h,\infty}\backslash \omega_h)\times\omega^{\tau}$ to the equation
\begin{equation}
 i\hbar s_N\bar{\partial}_t \Psi
 =\left(-c_{\hbar}{\Delta}_{hN}+V_{\infty}s_N\right)\overline{s}_t \Psi.
\label{n1}
\end{equation}
The boundary and initial conditions (\ref{l4}) and (\ref{l6}) imply that
\begin{equation}
 \left. \Psi^m\right|_{\Gamma_{h,\infty}\backslash \Gamma_h}=0\ \ \mbox{for}\ \ m\geq 1,\ \
 \Psi^0=0\ \ \mbox{on}\ \ \overline{\omega}_{h,\infty}\backslash (\omega_h\cup\Gamma_h).
\label{n2}
\end{equation}
\par Following \cite{AES03,DZ06,IZ11}, we apply the operator ${\mathcal F}_y$ to equation
(\ref{n1}) and obtain the 1D in space equation for the Fourier coefficients of $\Psi$
\begin{equation}
 i\hbar \Bigl(\sigma_q+\frac{h_x^2}{12}\partial_x\bar{\partial}_x\Bigr)\bar{\partial}_t\Psi^{(q)}
 =\Bigl[
 -c_{\hbar}(\sigma_q\partial_x\bar{\partial}_x
 -\lambda_q s_{Nx})
 +V_{\infty}\Bigl(\sigma_q+\frac{h_x^2}{12}\partial_x\bar{\partial}_x\Bigr)
 \Bigr]\overline{s}_t\Psi^{(q)}
\label{n5}
\end{equation}
on $\{x_j\}^{\infty}_{j=J} \times \omega^{\tau}$, with zero initial data
\[
 \Psi^{0(q)}=0\ \ \mbox{on}\ \ \{x_j\}^{\infty}_{j=J-1},
\]
see (\ref{n2}).
Dividing \eqref{n5} by $\sigma_q$, applying formulas \eqref{b2} and
\[
 \partial_x\bar{\partial}_x-\frac{\lambda_q}{\sigma_q}\Bigl(I+\frac{h_x^2}{12}\partial_x\bar{\partial}_x\Bigr)
 =\Bigl(I+\frac{1-\sigma_q}{\sigma_q^2}\lambda_q\frac{h_x^2}{12}\Bigr)\partial_x\bar{\partial}_x
 -\frac{\lambda_q}{\sigma_q}\Bigl(I+\frac{h_x^2}{12\sigma_q}\partial_x\bar{\partial}_x\Bigr)
\]
together with \eqref{b1}, we pass to the equation
\begin{equation}
 i\hbar s_{\theta_q x}\bar{\partial}_t\Psi^{(q)}
 =\left(-c_{\hbar,q}\partial_x\bar{\partial}_x
 +V_{\infty,q}s_{\theta_q x}\right)\overline{s}_t\Psi^{(q)}\ \
 \mbox{on}\ \ \{x_j\}^{\infty}_{j=J} \times \omega^{\tau},
\label{n5aa}
\end{equation}
with the coefficients $c_{\hbar,q}$ and $V_{\infty,q}$ and the parameter $\theta_q$ given by \eqref{l51}.
\par The discrete TBC for a 1D Schr\"odinger finite-difference equation with the average in $x$ like \eqref{n5aa} and constant coefficients was constructed in \cite{DZZ09}.
Taking into account that $\Psi^m\in H_h$ for any $m\geq 0$ and the Parseval equality, according to \cite{DZZ09} it can be represented as follows
\begin{equation}
\left.\left[
c_{\hbar,q}\bar{\partial}_x\overline{s}_t\Psi^{(q)}
-h_xs_{\theta_q x}^-
\bigl(i\hbar\,\bar{\partial}_t\Psi^{(q)}-V_{\infty,q}\overline{s}_t\Psi^{(q)}\bigr)
\right]\right|_{j=J}
=c_{\hbar,q}R_q *\Psi^{(q)}_J\ \ \mbox{on}\ \ \omega^{\tau},
\label{n5b}
\end{equation}
with the formulas for the kernel $R_q$ listed in the statement. (Actually the formulas are slightly modified and refined from misprints; also the
recently checked fixed sign in the formula for $c_{1q}$ when $0\leq\theta_q\leq\frac14$ is taken into account.)
\par Since $\Psi^{(q)}_j=\widetilde{\Psi}^{(q)}_j=\breve{\Psi}^{(q)}_j$ for $j=J-1$ and $J$, the discrete TBC \eqref{n5b} can be rewritten in the form \eqref{l5m} as well closely coupling it to the main equation \eqref{l2}.
\end{pf}
\par We need also to get an operator form of the derived discrete TBC.
\begin{proposition}
\label{pk2}
The discrete TBC \eqref{l5m} can be written in the operator form
\begin{gather}
 \bigl.{\mathcal D}_{xh}(\widetilde{\Psi}^m,\breve{\Psi}^m)\bigr|_{j=J}:=
 \Bigl.\Bigl[c_\hbar s_{Ny}\bar{\partial}_x\frac{\widetilde{\Psi}^{m}+\breve{\Psi}^{m}}{2}
 -h_xs_N^-\Bigl(i\hbar\,\frac{\widetilde{\Psi}^{m}-\breve{\Psi}^{m}}{\tau}
\Bigr.\Bigr.\Bigr.
\nonumber\\[1mm]
\Bigl.\Bigl.
-V_\infty\frac{\widetilde{\Psi}^{m}+\breve{\Psi}^{m}}{2}\Bigr)
 -h_xs_{Nx}^-c_\hbar{\partial}_y\bar{\partial}_y\frac{\widetilde{\Psi}^{m}+\breve{\Psi}^{m}}{2}
\Bigr]\Bigr|_{j=J}
=c_\hbar\mathcal{S}_{\rm ref}^m{\widetilde{\mathbf\Psi}}_J^m\ \ \text{for any}\ \ m\geq 1
\label{n5a}
\end{gather}
involving the operators
\begin{gather}
 s_{Nx}^{-}W_j:=\frac{1}{12}\,W_{j-1}+\frac{5}{12}\,W_j,\ \
 s_N^{-}:=s_{Nx}^{-}+\frac{h_y^2}{24}\partial_y\bar{\partial}_y,
\nonumber\\[1mm]
 \mathcal{S}_{\rm ref}^m{\mathbf\Psi}_J^m:=\mathcal{F}_y\bigl(\sigma_{N,q}R_q*\Psi^{(q)}_J\bigr),\ \
 \sigma_{N,q}:=\sigma_q\Bigl[1+\Bigl(\frac{h_xh_y\lambda_q}{12\sigma_q}\Bigr)^2\Bigr].
\label{n51}
\end{gather}
\end{proposition}
\begin{pf}
We go back to the discrete TBC \eqref{n5b} and perform the following transformations
\begin{gather*}
\sigma_q\bigl[c_{\hbar,q}\bar{\partial}_x\overline{s}_t\Psi^{(q)}
 -h_xs_{\theta_q x}^-\bigl(i\hbar\,\bar{\partial}_t\Psi^{(q)}-V_{\infty,q}\overline{s}_t\Psi^{(q)}\bigr)\bigr]
 =c_\hbar\Bigl[\sigma_q+\Bigl(\frac{1}{\sigma_q}-1\Bigr)\lambda_q\frac{h_x^2}{12}\Bigr]\bar{\partial}_x\overline{s}_t\Psi^{(q)}
\\[1mm]
 -h_x\Bigl(s_{Nx}^--\frac{h_y^2}{24}\lambda_q\Bigr)
  \bigl(i\hbar\,\bar{\partial}_t\Psi^{(q)}-V_\infty\overline{s}_t\Psi^{(q)}\bigr)
 +c_\hbar\lambda_q\Bigl(\frac{h_x}{2}-\frac{1}{\sigma_q}\frac{h_x^2}{12}\bar{\partial}_x\Bigr)\overline{s}_t\Psi^{(q)}
\\[1mm]
=c_\hbar\sigma_q\bar{\partial}_x\overline{s}_t\Psi^{(q)}
 -h_x\Bigl(s_{Nx}^--\frac{h_y^2}{24}\lambda_q\Bigr)
  \bigl(i\hbar\,\bar{\partial}_t\Psi^{(q)}-V_\infty\overline{s}_t\Psi^{(q)}\bigr)
 +h_xs_{Nx}^-c_\hbar\lambda_q\overline{s}_t\Psi^{(q)}
\end{gather*}
with the help of the formulas
\[
 \sigma_qs_{\theta_q x}^-=\frac{\sigma_q}{2}-\frac{h_x}{12}\bar{\partial}_x
 =\frac{1}{2}-\frac{h_x}{12}\bar{\partial}_x-\frac{h_y^2}{24}\lambda_q,\ \
 \frac{1}{2}-\frac{h_x}{12}\bar{\partial}_x=s_{Nx}^-.
\]
\par Thus applying $\mathcal{F}_y$ to \eqref{n5b} multiplied by $\sigma_q$, we obtain
\[
\Bigl.\Bigl[c_\hbar s_{Ny}\bar{\partial}_x\overline{s}_t\Psi
 -h_x\Bigl(s_{Nx}^-+\frac{h_y^2}{24}{\partial}_y\bar{\partial}_y\Bigr)
 \bigl(i\hbar\,\bar{\partial}_t\Psi-V_\infty\overline{s}_t\Psi\bigr)
 -h_xs_{Nx}^-c_\hbar{\partial}_y\bar{\partial}_y\overline{s}_t\Psi
\Bigr]\Bigr|_{j=J}
=c_\hbar\mathcal{S}_{\rm ref}{\mathbf\Psi}_J
\]
on $\omega^\tau$, where $\mathcal{S}_{\rm ref}{\mathbf\Psi}_J$ is given by \eqref{n51}.
As in the previous proof, this equation can be rewritten as \eqref{n5a}.
\end{pf}
\par The last form of the discrete TBC is in the spirit of \cite{DZ06}-\cite{DZZ13}, \cite{ZZ11,IZ11}
allowing to ensure both stability of schemes and the stable numerical implementation of the discrete TBCs.
\par Notice that the following important identity coupling the operators in the main equation (\ref{l2}) and the discrete TBC (\ref{n5a}) holds
\begin{gather}
\Bigl(
 i\hbar s_N\, \frac{\widetilde{\Psi}^m-\breve{\Psi}^m}{\tau}
 +c_{\hbar}{\Delta}_{hN}\frac{\widetilde{\Psi}^m+\breve{\Psi}^m}{2}
 -s_N\Bigl(\tilde{V}\frac{\widetilde{\Psi}^m+\breve{\Psi}^m}{2}\Bigr),W\Bigr)_{\omega_h}
\nonumber\\[1mm]
 -\Bigl(\bigl.
 {\mathcal D}_{xh}(\widetilde{\Psi}^m,\breve{\Psi}^m)
 \bigr|_{j=J},W_J
 \Bigr)_{\omega_{hy}}
 =i\hbar\Bigl(\frac{\widetilde{\Psi}^m-\breve{\Psi}^m}{\tau}
 -\tilde{V}\frac{\widetilde{\Psi}^m+\breve{\Psi}^m}{2},W \Bigr)_{\widetilde\omega_h,\,s_N}
\nonumber\\[1mm]
+c_{\hbar}
 \Bigl(
s_{Ny}\bar{\partial}_x\, \frac{\widetilde{\Psi}^m+\breve{\Psi}^m}{2}, \bar{\partial}_x W\Bigr)_{\tilde{\omega}_h}
+c_{\hbar}
 \Bigl(
-\partial_y\bar{\partial}_y\frac{\widetilde{\Psi}^m+\breve{\Psi}^m}{2},W\Bigr)_{\widetilde\omega_h,\,s_{Nx}},
\label{m1}
\end{gather}
for any $W$: $\overline{\omega}_h\rightarrow{\mathbb C}$ such that $W|_{j=0}=0$.
Here we have used the collection of $L^2$-mesh inner products
\begin{gather}
 (U,W)_{\omega_h}
 :=\sum_{j=1}^{J-1}\sum_{k=1}^{K-1}U_{j,k}W^*_{j,k}h_xh_y,\ \
 (U,W)_{\widetilde\omega_h}
 :=\sum_{j=1}^{J}\sum_{k=1}^{K-1}U_{j,k}W^*_{j,k}h_xh_y,
\nonumber\\
 (Z,\tilde{Z})_{\omega_{hy}}
 :=\sum_{k=1}^{K-1}Z_k\tilde{Z}^*_kh_y,\ \
 \left(U,W\right)_{\widetilde\omega_h,\,s}
 :=\left(sU,W\right)_{\omega_h}
 +\left(s^-U_J,W_J\right)_{\omega_{hy}}h_x,
\label{m2}
\end{gather}
with $s=s_N$ or $s_{Nx}$; notice that from \cite{DZZ09} it follows that the second sesquilinear form in (\ref{m2}) is Hermitian and positive definite on
functions $U,W$: $\overline{\omega}_h\to{\mathbb C}$ such that $U|_{\Gamma_h}=W|_{\Gamma_h}=0$.
We define the norms $\|\cdot\|_{\omega_h}$, $\|\cdot\|_{\tilde \omega_h}$ and $\|\cdot\|_{\omega_{hy}}$ associated to the first, second and third of these inner products.
\par Identity (\ref{m1}) appears after rearranging terms on the left-hand side and summing by parts with respect to $x$ in the term $c_{\hbar}s_{Ny}\partial_x\bar{\partial}_x$.

\begin{lemma}
\label{ln}
The operator ${\mathcal S}^m_{\rm ref}$ satisfies the inequality \cite{DZ06}
\begin{equation}
 \Ima\sum_{m=1}^M \left({\mathcal S}^m_{\rm ref}{\mathbf \Phi}^m,\overline{s}_t \Phi^m\right)_{\omega_{hy}}\tau \geq 0\ \
 \mbox{for any}\ \ M\geq 1,
\label{nn1}
\end{equation}
for any function $\Phi$: $\overline{\omega}_{hy}\times \overline{\omega}^{\,\tau}\to {\mathbb C}$ such that $\Phi^0=0$ and
$\Phi|_{k=0,K}=0$.
\end{lemma}
\begin{pf}
Following \cite{DZ06} and using the left formula \eqref{n51} and standard properties of ${\mathcal F}_y$, we obtain
\[
 \left({\mathcal S}^m_{\rm ref}{\mathbf \Phi}^m,\overline{s}_t \Phi^m\right)_{\omega_{hy}}
 =\frac{Y}{2}\sum_{q=1}^{K-1}\sigma_{N,q}
 \left(R_q*\Phi^{(q)} \right)^m\left(\overline{s}_t \Phi^m\right)^*.
\]
Consequently
\begin{equation}
 \Ima \sum_{m=1}^M \left({\mathcal S}^m_{\rm ref}{\mathbf \Phi}^m,\overline{s}_t \Phi^m\right)_{\omega_{hy}}\tau
 =\frac{Y}{2}
 \sum_{q=1}^{K-1}\sigma_{N,q}\Ima \sum_{m=1}^M
 \left(R_q*\Phi^{(q)}
 \right)^m
 \left(\overline{s}_t \Phi^m\right)^*\tau.
\label{nn1a}
\end{equation}
The result follows from the corresponding 1D inequality proved in \cite{DZZ09} (for any $\theta_q\leq\frac14$).
\end{pf}
By construction, the splitting scheme (\ref{l1})-(\ref{l5m}) on the finite space mesh has a solution. We need to prove its uniqueness. To this end, we assume that the auxiliary potential $\tilde{V}$ satisfies the condition
\begin{equation}
 |\tilde{V}(x')-\tilde{V}(x)|\leq L|x'-x|^\alpha\ \ \text{for any}\ \ 0\leq x\leq x'\leq X,
\label{condV}
\end{equation}
for some $\alpha\in [0,1]$.
It is simply the condition on the boundedness of spread in function values for $\alpha=0$,
the H\"{o}lder one for $\alpha\in (0,1)$ and the Lipschitz one for $\alpha=1$, on $[0,X]$.
Note that $L=0$ for $\tilde{V}=\textrm{const}$.
\begin{proposition}
\label{po}
Let $L\tau h_x^\alpha<8\hbar$.
Then the solution to the splitting scheme (\ref{l1})-(\ref{l5m}) on the finite space mesh is unique.
In addition, it satisfies the following $L^2$-stability bound
\begin{equation}
 \max_{0\leq m\leq M}\|\Psi^m\|_{\widetilde\omega_h}\leq \|\Psi^0_h\|_{\widetilde\omega_h}
 +\frac{6}{\hbar}\,\sum_{m=1}^M \| F^m\|_{\omega_h}\tau\ \ \ \text{for}\ \ m\geq 1
\label{o1}
\end{equation}
provided that $\Psi^0_h=0$ on $\Gamma_{1h}$ and
$\bigl.\Psi_{h}^0\bigr|_{j=J-1}=0$.
\end{proposition}
\begin{pf}
Assume that there exist two solutions of the scheme (\ref{l1})-(\ref{l5m}) and denote by $W$ their difference.
Obviously $W$ satisfies the homogeneous scheme (\ref{l1})-(\ref{l5m}), with $F=0$ and $\Psi^0_h=0$.
\par In order to establish uniqueness, it suffices to prove that if $W^0=0,\dots,W^{m-1}=0$, then $W^m=0$. Under this assumption $W^m$ satisfies a homogeneous equation
\begin{equation}
 i\frac{\hbar}{\tau}\, s_N W^m
 =-\frac{c_{\hbar}}{2}\,\Delta_{hN}W^m+\frac12s_N(\tilde{V}W^m)\ \ \mbox{on}\ \ \omega_h
\label{o2}
\end{equation}
supplemented with the homogeneous boundary conditions
\begin{equation}
 \left. W^m\right|_{\Gamma_h}=0,\ \ {\mathcal D}_{xh}\left(W^m,0\right)=c_{\hbar}{\mathcal S}^m_{\rm ref}{\mathbf W}^m_J\ \ \mbox{on}\ \ \Gamma_{1h},
\label{o3}
\end{equation}
where ${\mathbf W}^m=\{0,\dots,0,W^m\}$, with $0$ appearing $m$ times.
\par Applying the summation identity (\ref{m1}) in the case $\widetilde{\Psi}^m=W^m$, $\breve{\Psi}^m=0$ and $W=W^m$ as well as using (\ref{o2}) and (\ref{o3}), we get
\begin{gather}
 \Bigl(i\frac{\hbar}{\tau} W^m-\frac12\tilde{V}W^m,W^m\Bigr)_{\widetilde\omega_h,\,s_N}
 +\frac{c_{\hbar}}{2}
 \Bigl(s_{Ny}\, \bar{\partial}_xW^m,\bar{\partial}_xW^m\Bigr)_{\widetilde\omega_h}
\nonumber\\[1mm]
 +\frac{c_{\hbar}}{2}\left(-\partial_y\bar{\partial}_yW^m,W^m\right)_{\widetilde\omega_h,\,s_{Nx}}
 +c_{\hbar}
 \left(
 {\mathcal S}^m_{\rm ref}{\mathbf W}^m_J,W^m_J
 \right)_{\omega_{hy}}=0.
\label{q1}
\end{gather}
\par Let $H_{hy}$ be the space of functions $U$: $\overline{\omega}_{hy}\rightarrow{\mathbb C}$ such that
$U|_{k=0,K}=0$
endowed with the inner product $\left(\cdot,\cdot\right)_{\omega_{hy}}$.
Setting $(AP)|_{k=0,K}=0$
for $A=s_{Ny}$ and $-\partial_y\bar{\partial}_y$, we see that these operators are self-adjoint and positive definite in $H_{hy}$.
Therefore taking the imaginary part on the right of (\ref{q1}), we have
\[
 \frac{\hbar}{\tau}\left(W^m,W^m\right)_{\widetilde\omega_h,\,s_N}
 -\frac12\Ima \bigl(\tilde{V}W^m,W^m\bigr)_{\widetilde\omega_h,\,s_N}
 +c_{\hbar}\Ima \left({\mathcal S}^m_{\rm ref}{\mathbf W}^m_J,W^m_J\right)_{\omega_{hy}}=0.
\]
Due to independence of $\tilde{V}$ on $y$ and Lemma \ref{ln} we further derive
\[
 \frac{\hbar}{\tau}\left(W^m,W^m\right)_{\widetilde\omega_h,\,s_N}
 -\frac12\Ima \bigl(\tilde{V}W^m,W^m\bigr)_{\widetilde\omega_h,\,s_{Nx}}\leq 0.
\]
From \cite{DZZ09} it follows that
\[
 \frac13\|W^m\|_{\bar{\omega}_h}^2:=\frac13\Bigl(\|W^m\|_{\omega_h}^2+\frac{h_x}{2}\|W_J^m\|_{\omega_{hy}}^2\Bigr)
 \leq\left(W^m,W^m\right)_{\widetilde\omega_h,\,s_N}
\]
(compare with \eqref{h3a}) whereas according to Lemma 2.2 in \cite{ZL10} under condition \eqref{condV} we have
\[
 \bigl|\Ima\bigl(\tilde{V}W^m,W^m\bigr)_{\widetilde\omega_h,\,s_{Nx}}\bigr|\leq \frac{L}{12}h_x^\alpha\|W^m\|_{\bar{\omega}_h}^2.
\]
Consequently, $W^m=0$ provided that $L\tau h_x^\alpha<8\hbar$, and the uniqueness is proved.
\par Next, bound (\ref{o1}) clearly follows from the previous $L^2$-stability bound (\ref{l2}) in the case of the infinite space mesh since now $\Psi^0_h=0$ and $F^m=0$ on $\omega_{h,\infty}\backslash \omega_h$ for any $m\geq 1$.
\end{pf}
\par
Notice that, in the last proof, actually we have based upon a very particular case of inequality (\ref{nn1}), namely
\[
 0\leq\Ima\bigl({\mathcal F}_y^{-1}\bigl[\sigma_{N,q}R_q^0 U^{(q)}\bigr],U\bigr)_{\omega_{hy}}
 =\frac{Y}{2} \sum_{q=1}^{K-1}\sigma_{N,q}
 \Ima R_q^0 \bigl|U^{(q)}\bigr|^2,
\]
for any $U\in H_{hy}$ (see (\ref{nn1a})), which is clearly equivalent to
\[
 \Ima\, R_q^0\geq 0\ \ \mbox{for any}\ \ 1\leq q\leq K-1.
\]
\par The splitting scheme on the finite space mesh (\ref{l1})-(\ref{l5m}) can be effectively implemented (compare with \cite{DZZ13}).
Applying the operator $\sigma_q^{-1}{\mathcal F}_y$ to the main equation (\ref{l2}), similarly to the derivation of equation  \eqref{n5aa} we get a set of independent 1D finite-difference Schr\"odinger equations in $x$, for $1\leq q\leq K-1$
\begin{gather}
 i\hbar s_{\theta_qx}\,\frac{\widetilde{\Psi}^{m(q)}-\breve{\Psi}^{m(q)}}{\tau}
 =-c_{\hbar,q}\partial_x\bar{\partial}_x\frac{\widetilde{\Psi}^{m(q)}+\breve{\Psi}^{m(q)}}{2}
 +s_{\theta_qx}\Bigl(\tilde{V}_q
 \frac{\widetilde{\Psi}^{m(q)}+\breve{\Psi}^{m(q)}}{2}\Bigr)+\frac{F^{m(q)}}{\sigma_q}
\label{r1}
\end{gather}
on $\{x_j\}_{j=1}^{J-1}$, where $\tilde{V}_q:=\tilde{V}+c_{\hbar}\frac{\lambda_q}{\sigma_q}$.
This equation is supplemented with the boundary condition
\begin{gather}
 \widetilde{\Psi}^{m(q)}|_{j=0}=0
\label{r2}
\end{gather}
and the discrete TBC \eqref{l5m}.
\par Given $\Psi^{m-1}$, \textit{the direct algorithm} for computing $\Psi^{m}$ comprises five steps.
\begin{enumerate}
\item To compute explicitly $\breve{\Psi}^m={\mathcal E}^m\Psi^{m-1}$ on $\omega_h\cup\Gamma_{1h}$ (see (\ref{e6})).
\item To compute
$\breve{\Psi}^{m(q)}
=\bigl({\mathcal F}_y\breve{\Psi}^m\bigr)^{(q)}$ and
$F^{m(q)}=\left({\mathcal F}_yF^m\right)^{(q)}$,
for $1\leq q\leq K-1$.
\item  To compute
$\widetilde{\Psi}^{m(q)}$ by solving the 1D problems (\ref{r1}), (\ref{r2}) and \eqref{l5m}
for $1\leq q\leq K-1$; this includes the computation of the discrete convolutions on the right
of \eqref{l5m} so that $\widetilde{\Psi}_J^{1(q)},\dots,\widetilde{\Psi}_J^{m-1\,(q)}$ have to be stored.
\item  To compute
$\widetilde{\Psi}^{m}={\mathcal F}_y^{-1}\widetilde{\Psi}^{m(q)}$.
\item To compute explicitly
$\Psi^{m}={\mathcal E}^m\widetilde{\Psi}^{m}$ on $\omega_h\cup\Gamma_{1h}$ (see (\ref{e6})).
\end{enumerate}
\smallskip\par Steps 1 and 5 require $O\left(JK\right)$ arithmetic operations while
Steps 2 and 4 need $O\left(JK\log_2K\right)$ operations by using the fast Fourier transform (FFT) provided that
$K=2^{p}$, where $p$ is integer.
Step 3 requires $O((J+m)K)$ operations.
\par The total amount of arithmetic operations equals
$O\left(\left(J\log_2K+m\right)K\right)$ for computing the solution $\Psi^{m}$ at time level $m$ and
$O(\left(J\log_2K+M\right)KM)$ for computing one at all time levels $m=1,\dots,M$.
These amounts are the same as in \cite{DZZ13}.
\smallskip\par Notice that the given analysis is easily extended to the case of the problem in an infinite strip, with setting the following discrete TBC at the left boundary $x=0$ as well
\begin{gather}
\Bigl.\Bigl[
-c_{\hbar,q}{\partial}_x\,\frac{\widetilde{\Psi}^{m(q)}+\breve{\Psi}^{m(q)}}{2}
-h_xs_{\theta_q x}^+
\Bigl(i\hbar\,\frac{\widetilde{\Psi}^{m(q)}-\breve{\Psi}^{m(q)}}{\tau}
-V_{\infty,q}\,\frac{\breve{\Psi}^{m(q)}+\widetilde{\Psi}^{m(q)}}{2}\Bigr)
\Bigr]\Bigr|_{j=0}
\nonumber\\[1mm]
=c_{\hbar,q}\bigl(R_q *{\widetilde\Psi}^{(q)}_0\bigr)^m
\label{l5mm}
\end{gather}
for any $m\geq 1$ and $1\leq q\leq K-1$, where $s_{\theta x}^+W_j:=(\frac12-\theta) W_j+\theta W_{j+1}$ (for simplicity, we suppose that $V(x,y)=V_\infty$ for $x\leq h_x$ though clearly $V_{\pm\infty}$ could be different).
\section{\large Numerical experiments}
\setcounter{equation}{0}
\setcounter{proposition}{0}
\setcounter{theorem}{0}
\setcounter{lemma}{0}
\setcounter{corollary}{0}
\setcounter{remark}{0}
The above described direct algorithm has been implemented. For numerical experiments, we take $\hbar=1$ and $c_\hbar=1$.
We solve the initial-boundary value problem in the infinite strip $\mathbb{R}\times (0,Y)$ taking the computational domain $\bar{\Omega}_{X,Y}\times[0,T]$, where $\Omega_{X,Y}=(0,X)\times (0,Y)$.
\par Let the initial function be the standard Gaussian wave package
\[
 \psi^0(x,y)=
 \psi_G(x,y) :=
 \exp\left\{ ik(x-x^{(0)})-\frac{(x-x^{(0)})^2+(y-y^{(0)})^2}
 {4\alpha}\right\}\ \ \text{on}\ \ \mathbb{R}^2.
\]
We take its parameters $k=30\sqrt{2}$ (the wave number), $\alpha=\frac{1}{120}$ and $(x^{(0)},y^{(0)})=(1,\frac{Y}{2})$ (the position of the modulus' maximum) like in \cite{DZZ13}.
\par We respectively modify the splitting in potential scheme  (\ref{l1})-(\ref{l5m}) enlarging $\omega_h\cup \Gamma_{1h}$ by $\{(0,y_k);\, 1\leq k\leq K-1\}$ in \eqref{l1} and \eqref{l3} and replacing the boundary conditions \eqref{l4} by
\[
 \breve{\Psi}^m|_{k=0,K}=0,\ \ \widetilde{\Psi}^m|_{k=0,K}=0,\ \
 \Psi^m|_{k=0,K}=0
\]
together with the left discrete TBC \eqref{l5mm}.
Recall that we use the uniform meshes in $x$, $y$ and $t$
with the steps correspondingly $h_x=\frac{X}{J}$, $h_y=\frac{Y}{K}$ and $\tau=\frac{T}{M}$.
\smallskip\par\textbf{Example A.}
We first take a modified P\"{o}schl-Teller \cite{F71} potential (a barrier)
\[
 V(x)=\frac{\alpha_0^2 c_1}{\cosh^2 \alpha_0(x-x^*)}
\]
depending only on $x$.
We set $\alpha_0=6$, $c_1=47$ and $x^*=2$. Then $\max_\mathbb{R}V(x)=V(x^*)=\alpha_0^2 c_1=1692$, and though the potential is smooth, its derivatives in $x$ are rather large.

\par We choose $(X,Y)=(4,4.2)$, then both $V$ and $\psi_G$ are small enough,
namely, $|V(x)|<2.6\cdot10^{-7}$ and $|\psi_G(x,y)|<9.4\cdot 10^{-14}$ outside $\bar{\Omega}_{X,Y}$.
Let also $T=t_M=0.05$.
\par
\begin{figure}[ht]
\begin{multicols}{2}
    \includegraphics[scale=0.5]{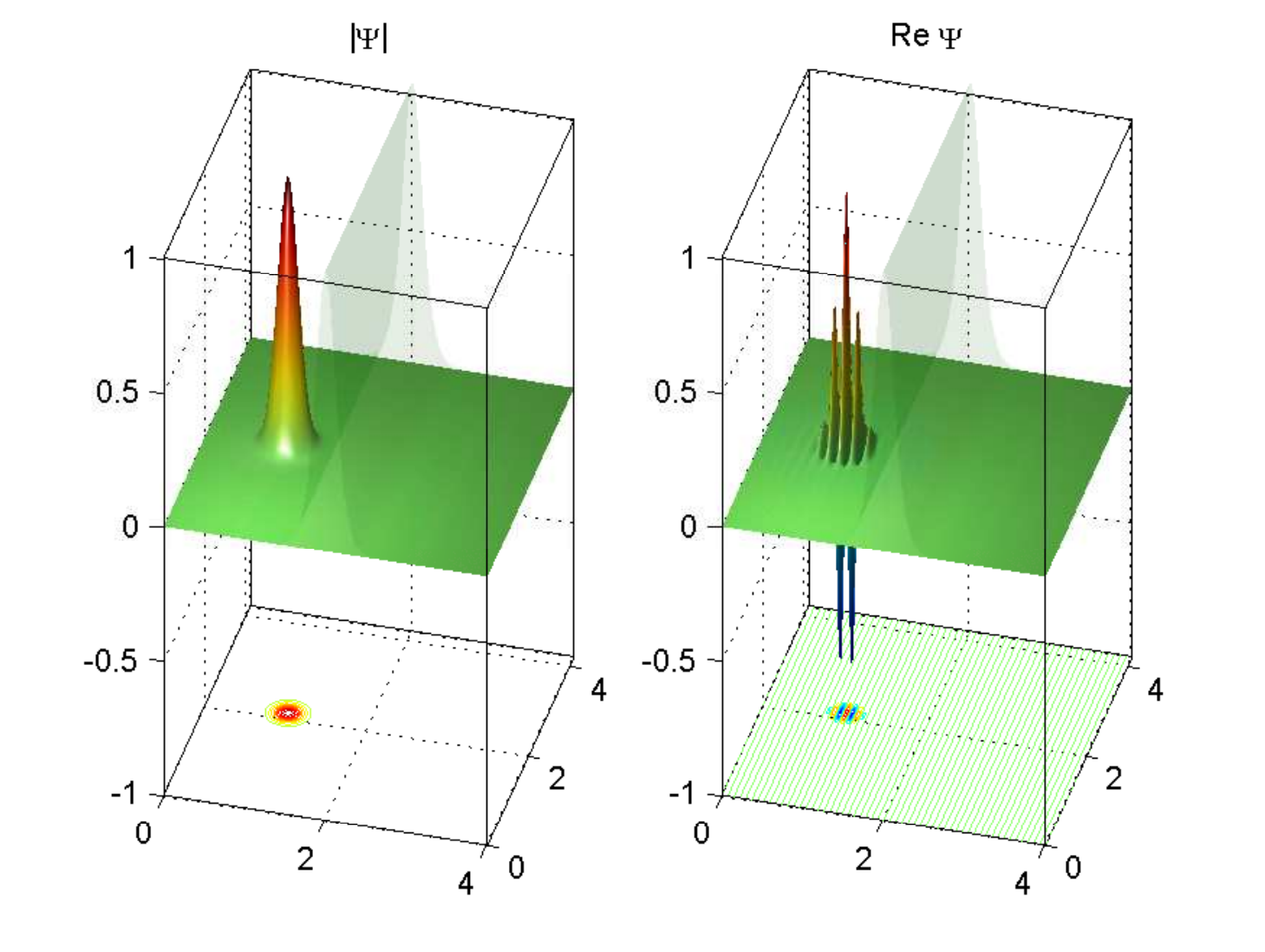}\vspace{0cm}\\
   \centerline{\small{$m=0$}}\\
    \includegraphics[scale=0.5]{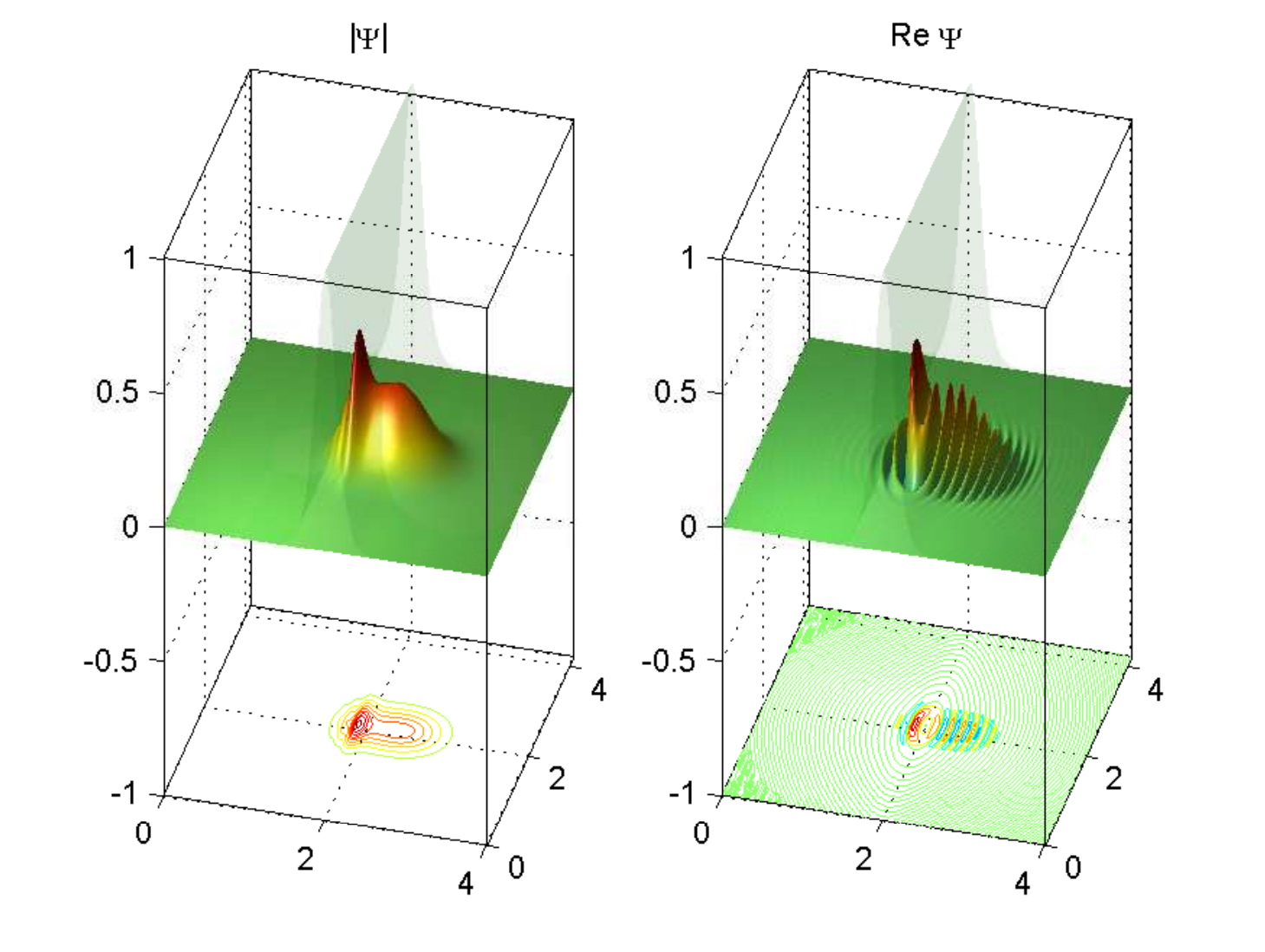}\vspace{0cm}\\
     \centerline{\small{$m=360$}}\\
\end{multicols}
\begin{multicols}{2}
    \includegraphics[scale=0.5]{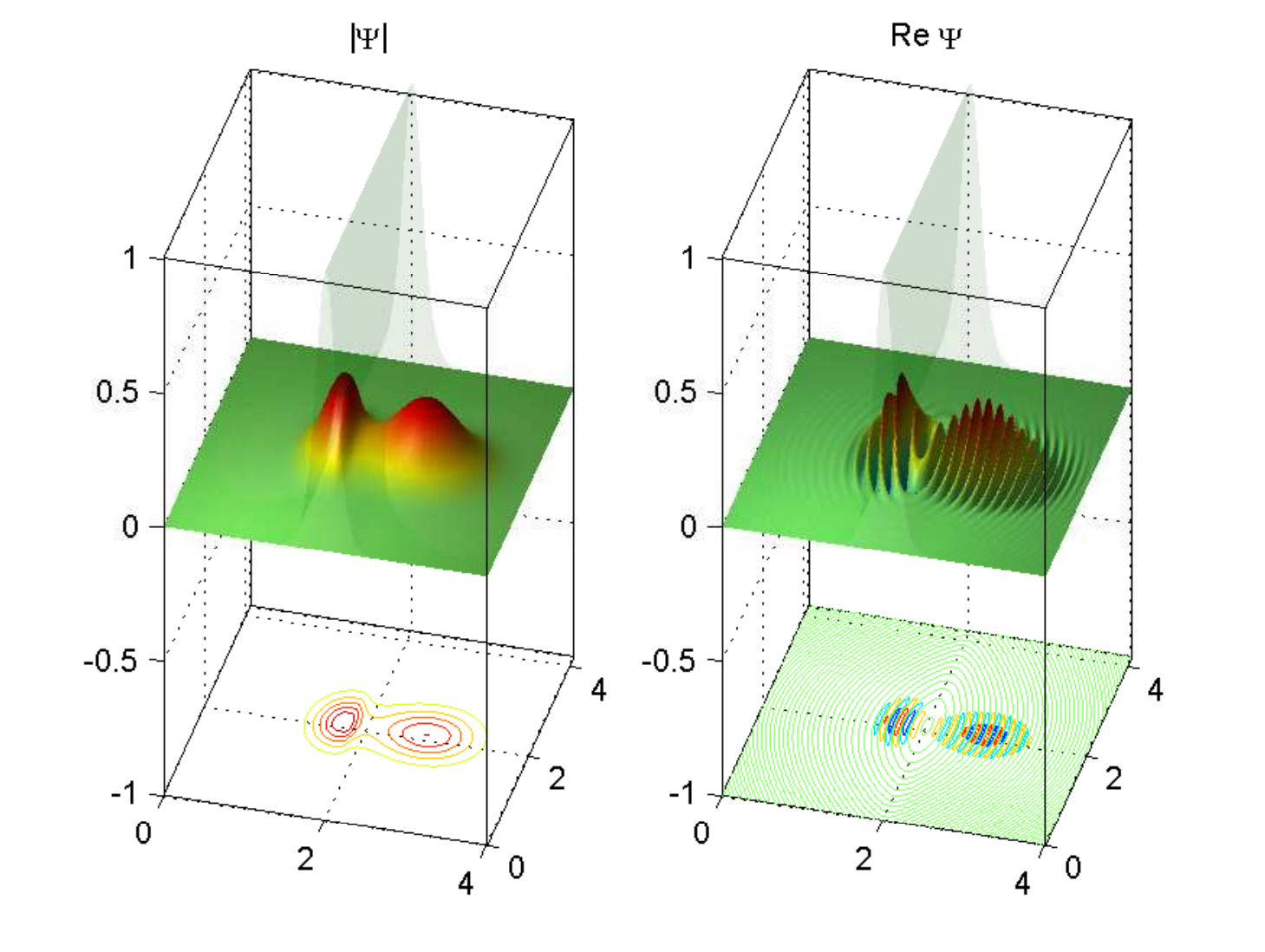}\vspace{0cm}\\%
    \centerline{\small{$m=440$}}\\
    \includegraphics[scale=0.5]{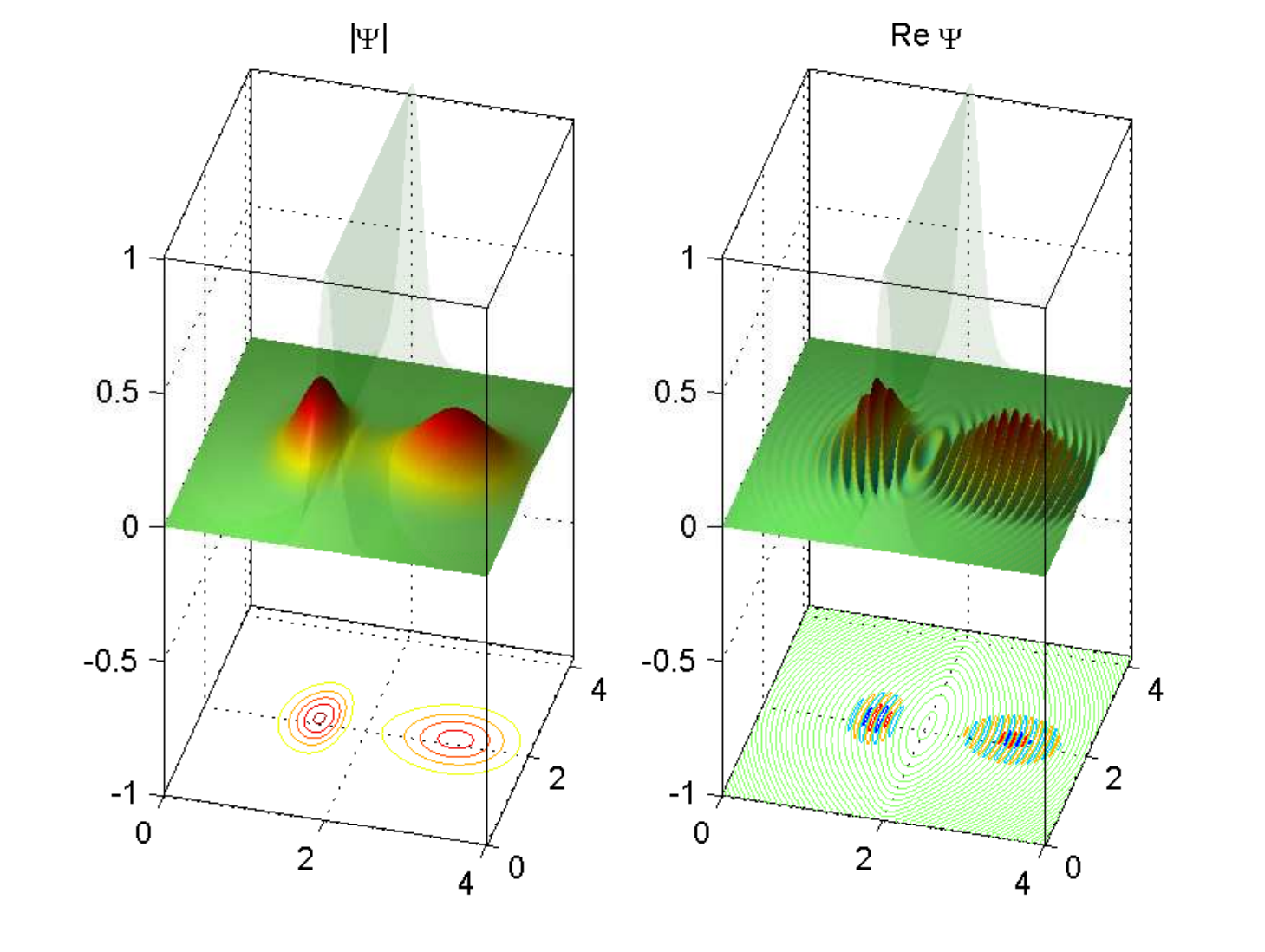}\vspace{0cm}\\%
     \centerline{\small{$m=520$}}\\
\end{multicols}
\caption{\small{Example A. The modulus and the real part of the numerical solution $\Psi^m$, $m=0, 360, 440,$ and $520$, for $(J,K,M)=(400,64,1000)$}}
\label{SSP:EX22b:B:Solution1}
\end{figure}

 The modulus and the real part of $\psi_G$ together with the normalized barrier (when $\alpha_0^2c_1=1$) are shown on the computational domain on Figure \ref{SSP:EX22b:B:Solution1}, for $m=0$.
\par
\begin{figure}[ht]
\begin{multicols}{2}
    \includegraphics[scale=0.5]{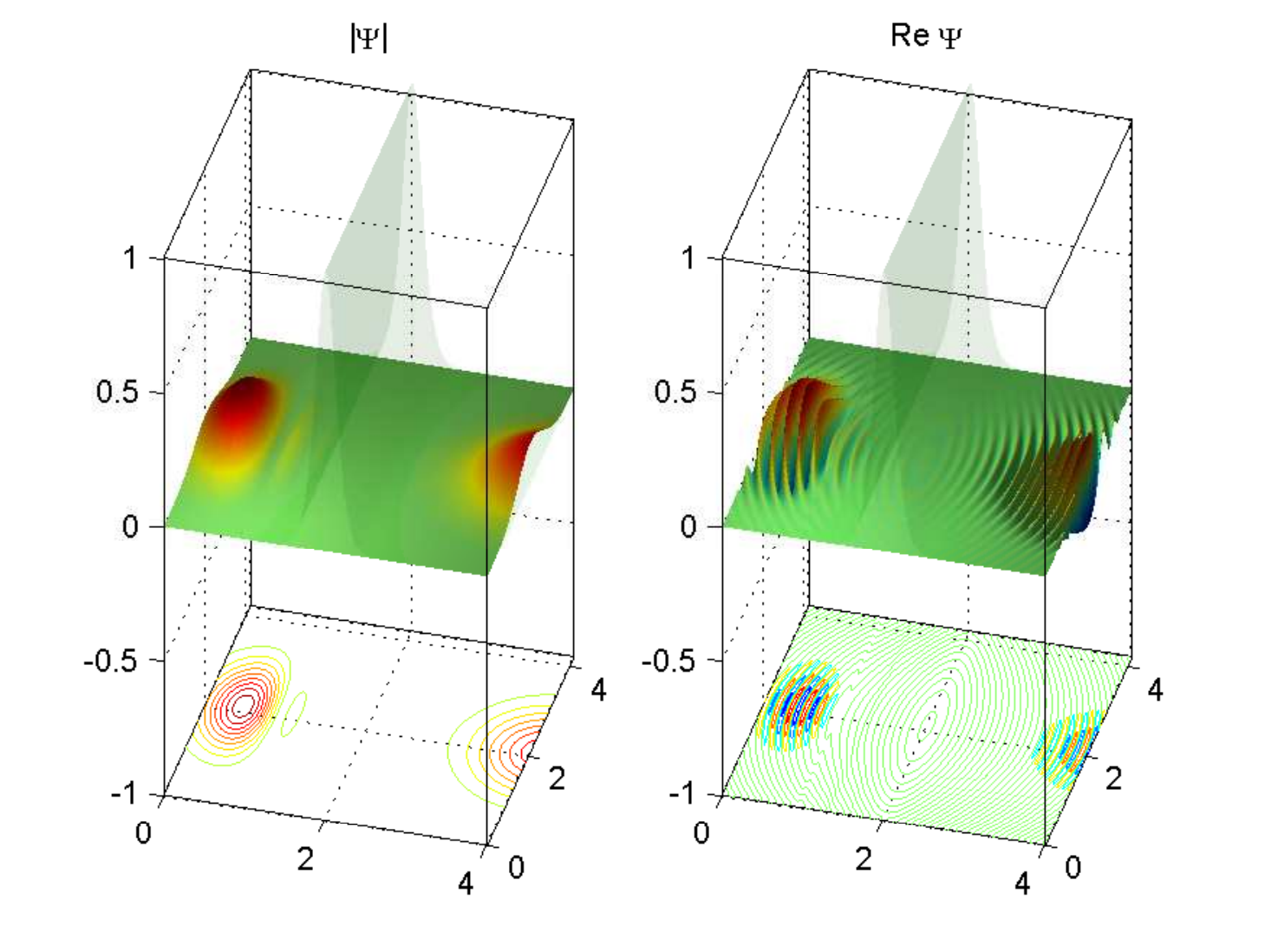}\vspace{0cm}\\%
     \centerline{\small{$m=780$}}\\
    \includegraphics[scale=0.5]{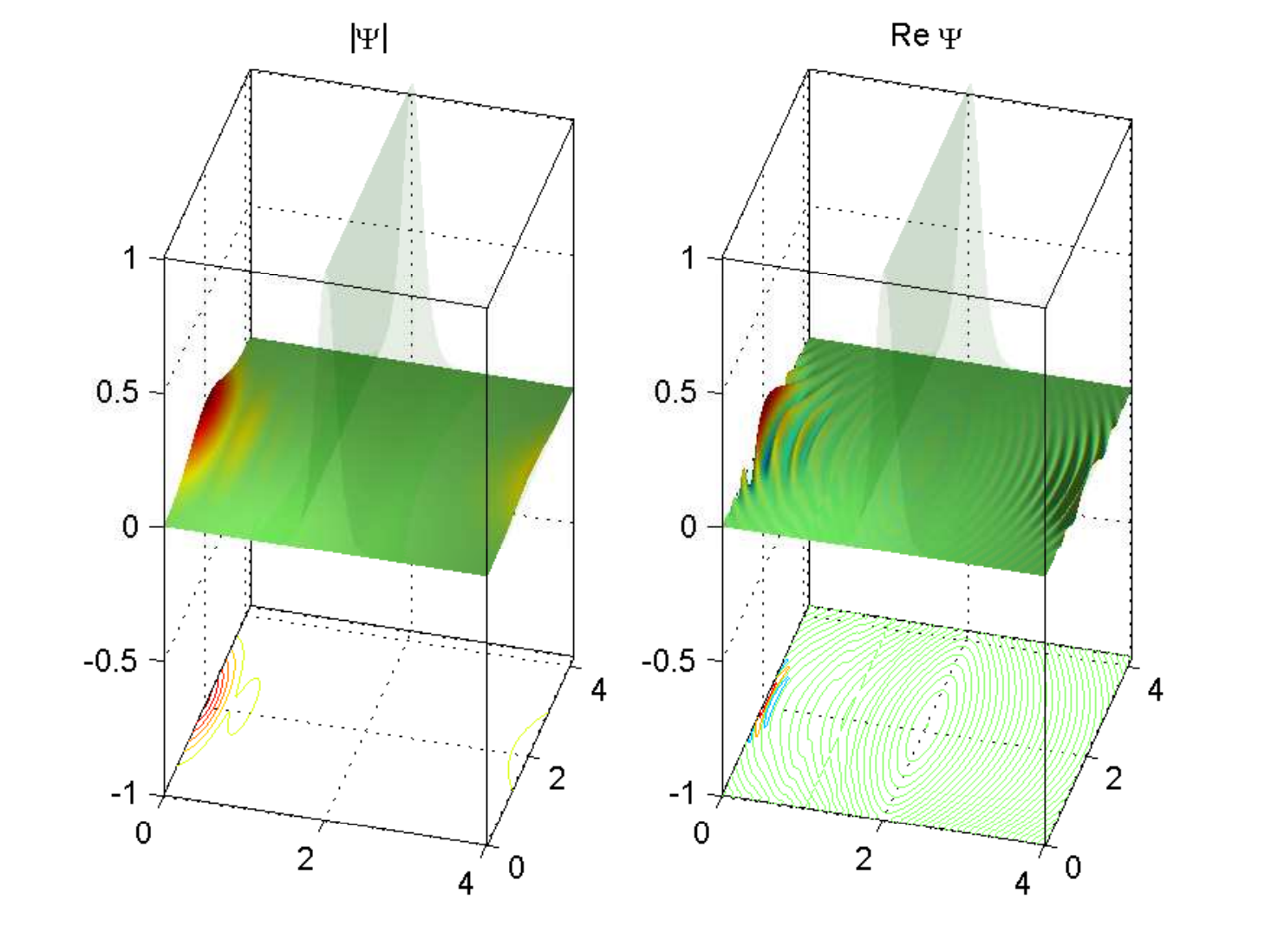}\vspace{0cm}\\%
     \centerline{\small{$m=960$}}\\
\end{multicols}
\caption{\small{Example A. The modulus and the real part of the numerical solution $\Psi^m$, $m=780$ and $960,$ for $(J,K,M)=(400,64,1000)$}
\label{SSP:EX22b:B:Solution2}}
\end{figure}

 We put $\tilde{V}=0$ and $\Delta V=V$ and compute the numerical solution $\Psi^{m}$ for $(J,K,M)=(400,64,1000)$
so that $h_x=10^{-2}$, $h_y\approx6.56\cdot10^{-2}$ and $\tau=5\cdot 10^{-5}$.
Its modulus and real part are presented on Figures \ref{SSP:EX22b:B:Solution1} and \ref{SSP:EX22b:B:Solution2}, for the time moments $t_m=m\tau$, $m=360, 440, 520, 780$ and $960$.
The wave package is separated by the barrier into two comparable reflected and transmitted parts moving in opposite $x$-directions and leaving the computational domain. Notice the oscillating behavior of the real part of the wave.
\par
\begin{figure}[ht]\centerline{
   \includegraphics[width=0.5\linewidth]{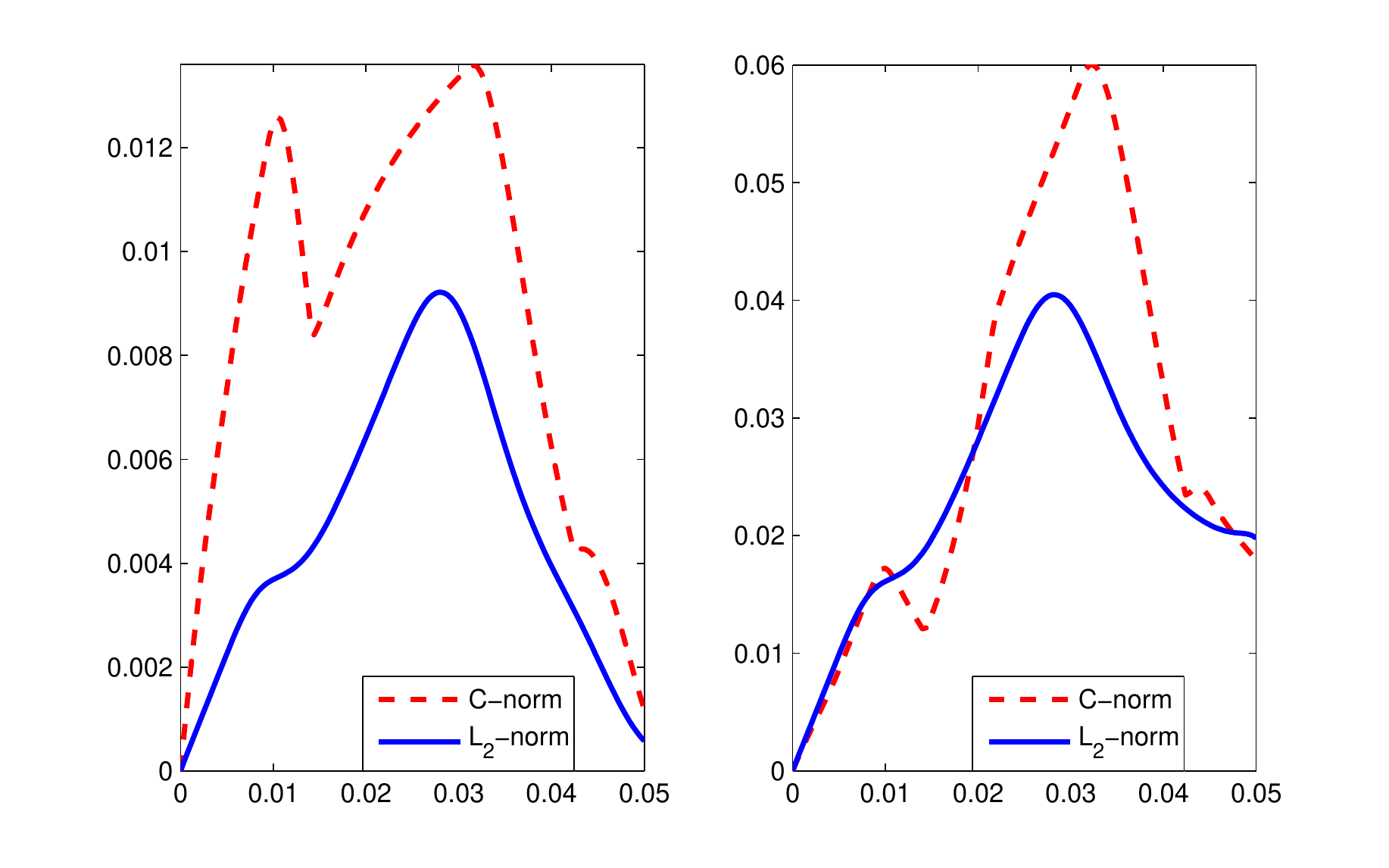}}
\caption{\small{Example A. The absolute and relative differences in $C$ and $L^2$ norms between the numerical solutions for $(J,K,M)=(400,64,1000)$ and $(1600,256,4000)$, in dependence with time}
\label{SSP:EX22b:B:Errors}}
\end{figure}

 To justify the choice of the mesh with not so large $J$ and $K$, on Figure \ref{SSP:EX22b:B:Errors} we
give the absolute and relative differences in $C$ and $L^2$ space mesh norms for the numerical solutions for $(J,K,M)=(400,64,1000)$ and $(1600,256,4000)$ that all three are 4 times larger, in dependence with time. The level of differences is enough to represent the correct graphs of the exact solution.

\par We intend to study the error behavior in more detail. Given $(J,K,M)$, denote by $\Psi_{J,K,M}$ the corresponding numerical solution. Let the following error representation
\[
 \psi-\Psi_{J,K,M}=ah_x^\alpha+bh_y^\beta+c\tau^\gamma+\delta(h_x,h_y,\tau)
\]
with the first main error term and the residual $\delta$ such that
\[
 \|\delta(h_x,h_y,\tau)\|=o(h_x^\alpha+h_y^\beta+\tau^\gamma)\ \ \text{as}\ \ h_x+h_y+\tau\to 0
\]
be valid, for some positive exponents $\alpha$, $\beta$ and $\gamma$, where the functions $a$, $b$ and $c$ are non-zero and mesh independent as well as $\|\cdot\|$ is a space-time mesh norm. Given $(\bar{J},\bar{K},\bar{M})$ and corresponding $\bar{h}_x=\frac{X}{\bar{J}}$, $\bar{h}_y=\frac{Y}{\bar{K}}$ and $\bar{\tau}=\frac{T}{\bar{M}}$ together with integer $\ell\geq 1$, derive
\[
 \rho_{x,\bar{J},\bar{K},\bar{M},2^\ell}:=\|\Psi_{\bar{J},\bar{K},\bar{M}}-\Psi_{\bar{J}/2^\ell,\bar{K},\bar{M}}\|
 =\|a\|(2^\ell-1)\bar{h}_x^\alpha+o(\bar{h}_x^\alpha+\bar{h}_y^\beta+\bar{\tau}^\gamma).
\]
Let $\bar{J}$, $\bar{K}$, $\bar{M}$ be sufficiently large. Therefore we get
\begin{equation}
 R_{x,\bar{J},\bar{K},\bar{M},2^\ell}:=\frac{\rho_{x,\bar{J},\bar{K},\bar{M},2^{\ell+1}}}{\rho_{x,\bar{J},\bar{K},\bar{M},2^\ell}}
 \sim r_{\alpha,\ell}:=\frac{2^{(\ell+1)\alpha}-1}{2^{\ell\alpha}-1}
\label{num3}
\end{equation}
provided that $\bar{h}_y^\beta=O(\bar{h}_x^\alpha)$ and $\bar{\tau}^\gamma=O(\bar{h}_x^\alpha)$.
Similarly
\[
 \rho_{y,\bar{J},\bar{K},\bar{M},2^\ell}:=\|\Psi_{\bar{J},\bar{K},\bar{M}}-\Psi_{\bar{J},\bar{K}/2^\ell,\bar{M}}\|
 =\|b\|(2^\ell-1)\bar{h}_y^\beta+o(\bar{h}_x^\alpha+\bar{h}_y^\beta+\bar{\tau}^\gamma)
\]
and
\begin{equation}
 R_{y,\bar{J},\bar{K},\bar{M},2^\ell}:=\frac{\rho_{y,\bar{J},\bar{K},\bar{M},2^{\ell+1}}}{\rho_{y,\bar{J},\bar{K},\bar{M},2^\ell}}
 \sim r_{\beta,\ell}
\label{num5}
\end{equation}
provided that $\bar{h}_x^\alpha=O(\bar{h}_y^\beta)$ and $\bar{\tau}^\gamma=O(\bar{h}_y^\beta)$ as well as
\[
 \rho_{t,\bar{J},\bar{K},\bar{M},2^\ell}:=\|\Psi_{\bar{J},\bar{K},\bar{M}}-\Psi_{\bar{J},\bar{K},\bar{M}/2^\ell}\|
 =\|c\|(2^\ell-1)\bar{\tau}^\gamma+o(\bar{h}_x^\alpha+\bar{h}_y^\beta+\bar{\tau}^\gamma)
\]
and
\begin{equation}
 R_{t,\bar{J},\bar{K},\bar{M},2^\ell}:=\frac{\rho_{t,\bar{J},\bar{K},\bar{M},2^{\ell+1}}}{\rho_{t,\bar{J},\bar{K},\bar{M},2^\ell}}
 \sim r_{\gamma,\ell}
\label{num7}
\end{equation}
provided that $\bar{h}_x^\alpha=O(\bar{\tau}^\gamma)$ and $\bar{h}_y^\beta=O(\bar{\tau}^\gamma)$.
The quantities $R_{x,\bar{J},\bar{K},\bar{M},2^\ell}$, $R_{y,\bar{J},\bar{K},\bar{M},2^\ell}$ and $R_{t,\bar{J},\bar{K},\bar{M},2^\ell}$ are computable, and by analyzing their behavior one can try to judge on the practical values of
the exponents  $\alpha$, $\beta$ and $\gamma$ when the exact solution is unknown (that is quite standard).
\par In this respect, notice that
\begin{gather}
 r_{\alpha,4}\approx 4.01,\ \ r_{\alpha,3}\approx 4.05,\ \ r_{\alpha,2}=4.2,\ \ r_{\alpha,1}=5\ \ \text{for}\ \ \alpha=2,
\label{num13}\\[1mm]
 r_{\alpha,4}\approx 8.00,\ \ r_{\alpha,3}\approx 8.01,\ \ r_{\alpha,2}\approx 8.11\ \ r_{\alpha,1}=9\ \ \text{for}\ \ \alpha=3,
\label{num14}\\[1mm]
 r_{\alpha,4}\approx 16.00,\ \ r_{\alpha,3}\approx 16.00,\ \ r_{\alpha,2}=16.06,\ \ r_{\alpha,1}=17\ \ \text{for}\ \ \alpha=4.
\label{num15}
\end{gather}
\par We implement this approach. In Table \ref{SSP:tab:EX22b:M}, the quantities
\begin{gather*}
 \rho_{x,\bar{J},\bar{K},\bar{M},2^\ell},\ \ R_{x,\bar{J},\bar{K},\bar{M},2^\ell} \ \
 \text{for}\ \ (\bar{J},\bar{K},\bar{M})=(3200,256,4444)\ \ \text{and}\ \ \ell=4,3,2,1,\\[1mm]
 \rho_{y,\bar{J},\bar{K},\bar{M},2^\ell},\ \ R_{y,\bar{J},\bar{K},\bar{M},2^\ell} \ \
 \text{for}\ \ (\bar{J},\bar{K},\bar{M})=(1600,512,4444)\ \ \text{and}\ \ \ell=4,3,2,1,\\[1mm] \rho_{t,\bar{J},\bar{K},\bar{M},2^\ell},\ \ R_{t,\bar{J},\bar{K},\bar{M},2^\ell} \ \
 \text{for}\ \ (\bar{J},\bar{K},\bar{M})=(1600,256,4000)\ \ \text{and}\ \ \ell=4,3,2,1
\end{gather*}
are sequentially presented both for the mesh norm $\|\cdot\|=E_C$ uniform in time and space and $\|\cdot\|=E_{L^2}$ uniform in time and $L^2$ in space. Comparing the data with \eqref{num13}-\eqref{num15}, we conclude on $\alpha\approx 4$, $3<\beta\lesssim 4$ and $\gamma\approx 2$ that is in agreement with their expected theoretical values.
On the other hand, this is not so clear a priori since the derivatives of both the oscillating exact solution and the potential are rather large.
 Notice also the low error values in $J$ and $K$ for $\ell=2,1$, and the most close correspondence between $R_C$ for $M$ and \eqref{num13}.
\begin{table}
\begin{center}
\begin{tabular}{|r|c|c|c|c|c|}
  \hline
  $J$ & $E_C$  & $E_{L^2}$  & $R_C$     & $R_{L^2}$ \\
  \hline
  200  & $0.411$$E$$-1$  &$0.208$$E$$-1$    & -        & -      \\
  400  & $0.240$$E$$-2$  &$0.120$$E$$-2$    & 17.13    & 17.33  \\
  800  & $0.147$$E$$-3$  &$0.664$$E$$-4$    & 16.32    & 18.07  \\
  1600 & $0.897$$E$$-5$  &$0.564$$E$$-5$    & 16.39    & 11.77  \\

  \hline
\end{tabular}
\end{center}

\begin{center}
\begin{tabular}{|r|c|c|c|c|c|}
  \hline
   $K$ & $E_C$  & $E_{L^2}$  & $R_C$     & $R_{L^2}$ \\
  \hline
  32  & $0.227$$E$$-1$  &$0.191$$E$$-1$    & -        & -                \\
  64  & $0.220$$E$$-2$  &$0.150$$E$$-2$    & 10.32    & 12.73            \\
  128 & $0.129$$E$$-3$ & $0.785$$E$$-4$    & 17.05    & 19.11            \\
  256 & $0.993$$E$$-5$ & $0.789$$E$$-5$    & 12.99    & \hphantom{1}9.95 \\
  \hline
\end{tabular}
\end{center}

\begin{center}
\begin{tabular}{|r|c|c|c|c|c|}
  \hline
   $M$ & $E_C$  & $E_{L^2}$  & $R_C$     & $R_{L^2}$ \\
  \hline
  250  & $0.206$\hphantom{$E$$-1$}     &$0.141$\hphantom{$E$$-1$}       &  -      &  -       \\
  500  & $0.508$$E$$-1$                &$0.371$$E$$-1$                  & 4.05    & 3.80     \\
  1000 & $0.121$$E$$-1$                &$0.880$$E$$-2$                  & 4.21    & 4.22     \\
  2000 & $0.240$$E$$-2$                &$0.180$$E$$-2$                  & 5.04    & 4.89     \\
  \hline
\end{tabular}
\end{center}
\caption{Example A. The change in numerical solution in maximum in time $C$ and $L^2$ space norms for redoubling
$J$ (and the reference $(\bar{J},\bar{K},\bar{M})=(3200,256,4444)$),
or $K$ (and the reference $(\bar{J},\bar{K},\bar{M})=(1600,512,4444)$),
or $M$ (and the reference $(\bar{J},\bar{K},\bar{M})=(1600,256,4000)$)}
\label{SSP:tab:EX22b:M}
\end{table}
\par
\begin{figure}[ht]
\begin{multicols}{2}
    \includegraphics[scale=0.37]{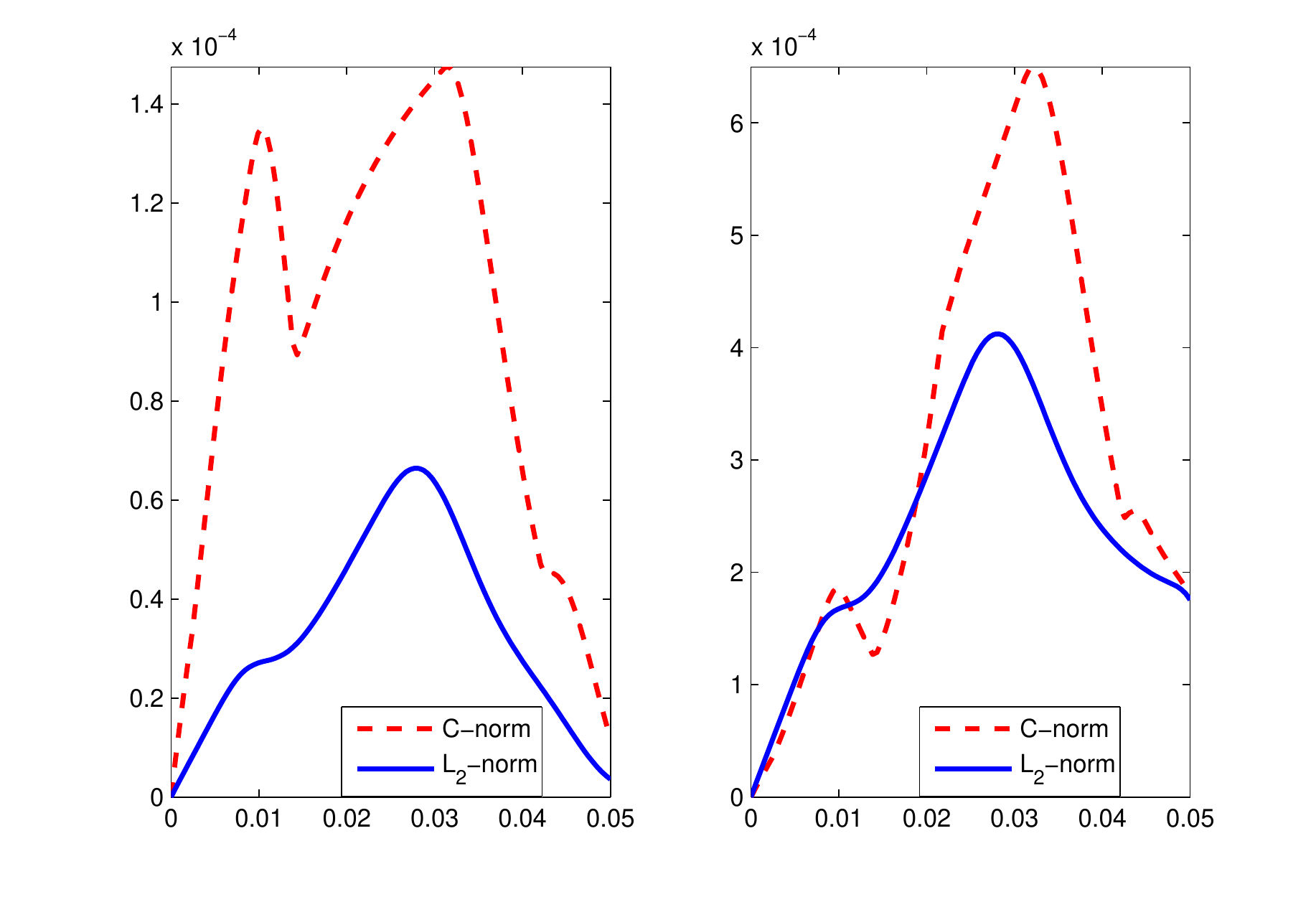}\vspace{0cm}\\
\caption{\small{Example A.
The absolute and relative differences in $C$ and $L^2$ norms between the numerical solutions for
$J=800$ and 3200, while $(K,M)=(256,4444)$, in dependence with time}}
\label{SSP:EX22b:B:Solution_J}
    \includegraphics[scale=0.42]{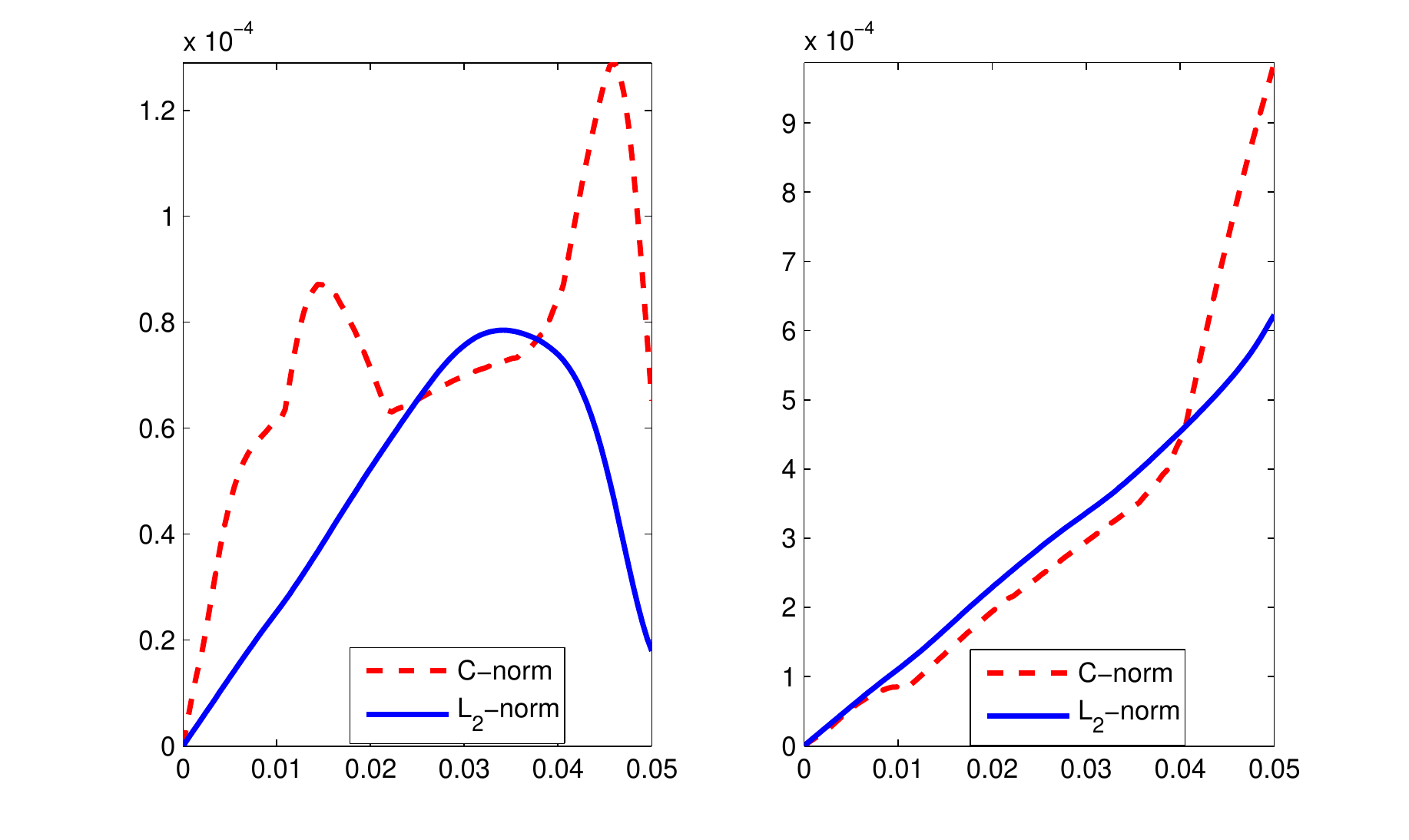}\vspace{0cm}\\
\caption{\small{Example A.
The absolute and relative differences in $C$ and $L^2$ norms between the numerical solutions for
$K=128$ and 512, while $(J,M)=(1600,4444)$, in dependence with time}}
\label{SSP:EX22b:B:Solution_K}
\end{multicols}
\end{figure}

 In addition, Figures \ref{SSP:EX22b:B:Solution_J} and \ref{SSP:EX22b:B:Solution_K} present the typical behavior of $\Psi_{\bar{J},\bar{K},\bar{M}}-\Psi_{\bar{J}/2^\ell,\bar{K},\bar{M}}$ and $\Psi_{\bar{J},\bar{K},\bar{M}}-\Psi_{\bar{J},\bar{K}/2^\ell,\bar{M}}$
both in $C$ and $L^2$ space norms, in absolute and relative forms, in dependence with time, for above $(\bar{J},\bar{K},\bar{M})$ and $\ell=2$.

\begin{figure}[ht]
\begin{multicols}{2}
    \includegraphics[scale=0.5]{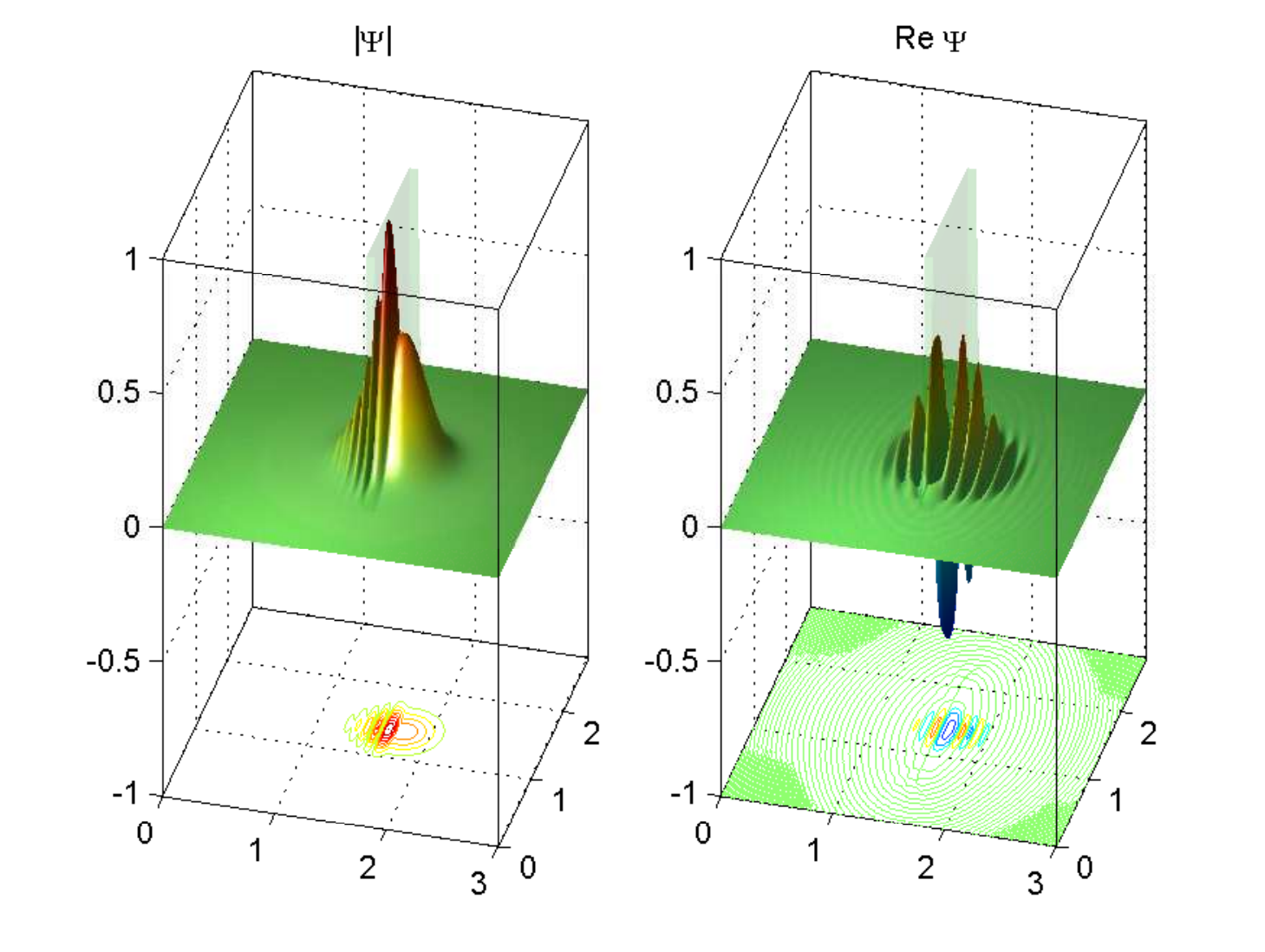}\vspace{0cm}\\%
\centerline{\small{$m=204$}}\\
    \includegraphics[scale=0.5]{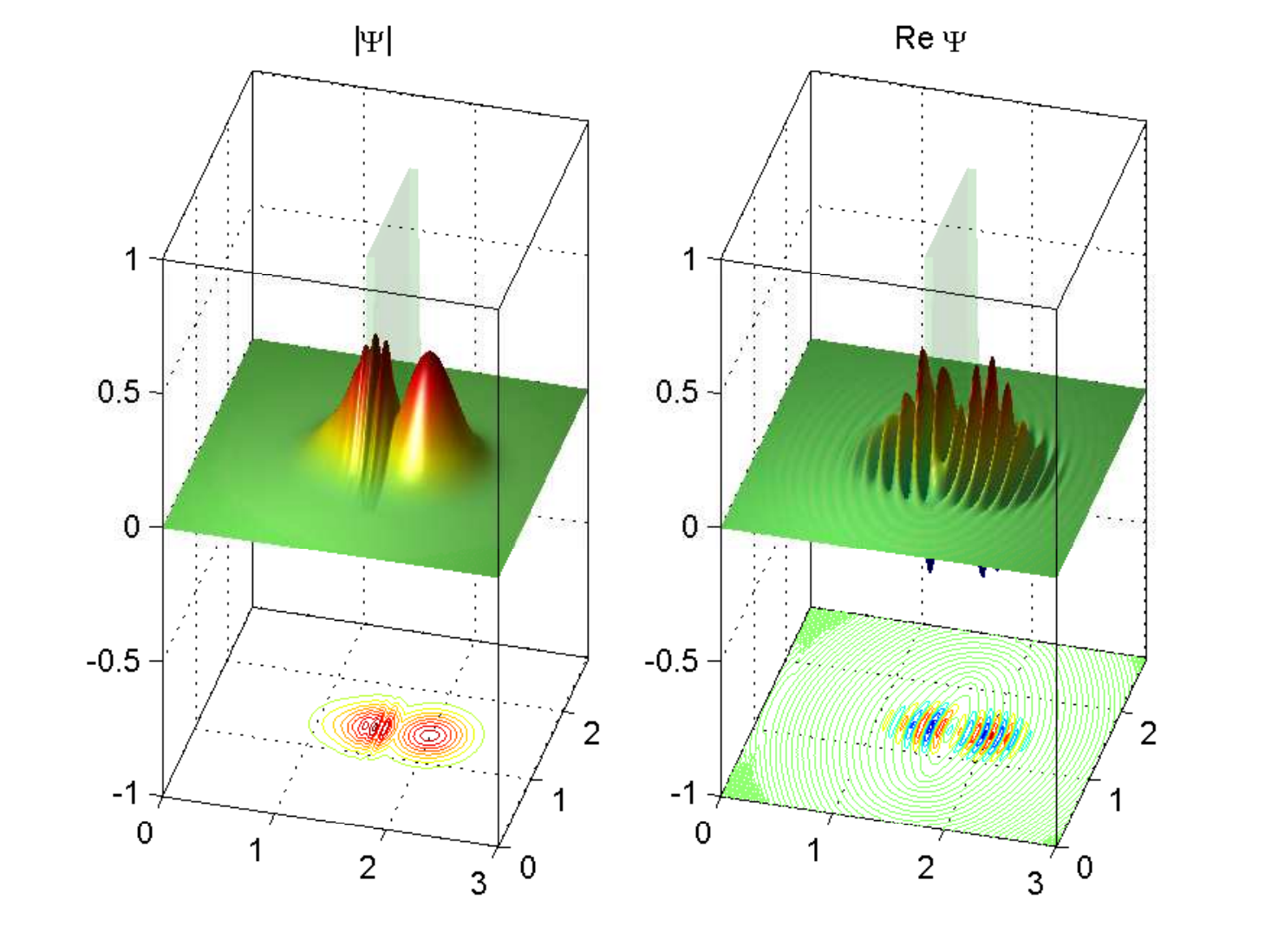}\vspace{0cm}\\%
\centerline{\small{$m=264$}}\\
\end{multicols}
    \begin{multicols}{2}
    \includegraphics[scale=0.5]{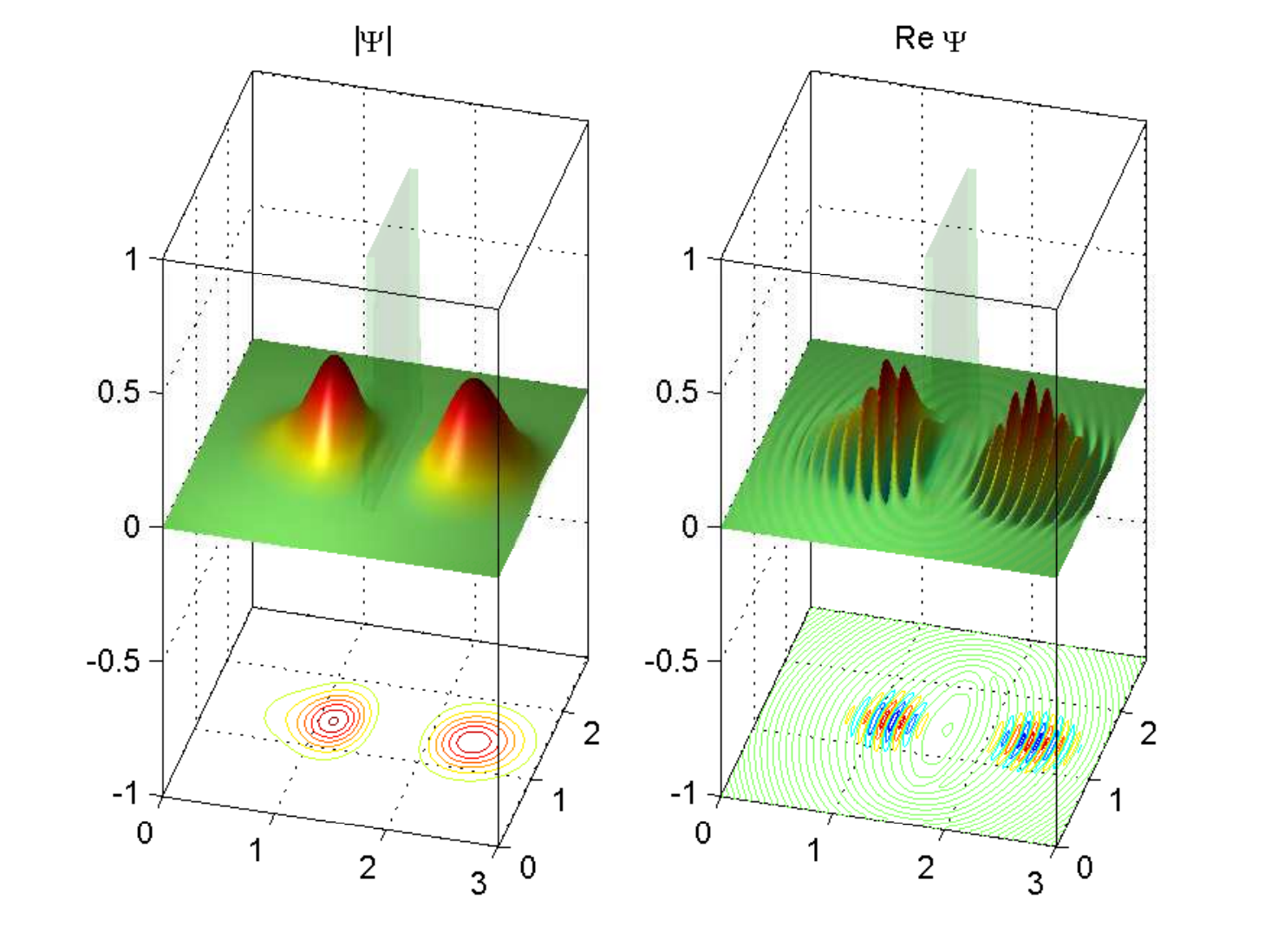}\vspace{0cm}\\%
\centerline{\small{$m=360$}}\\
    \includegraphics[scale=0.5]{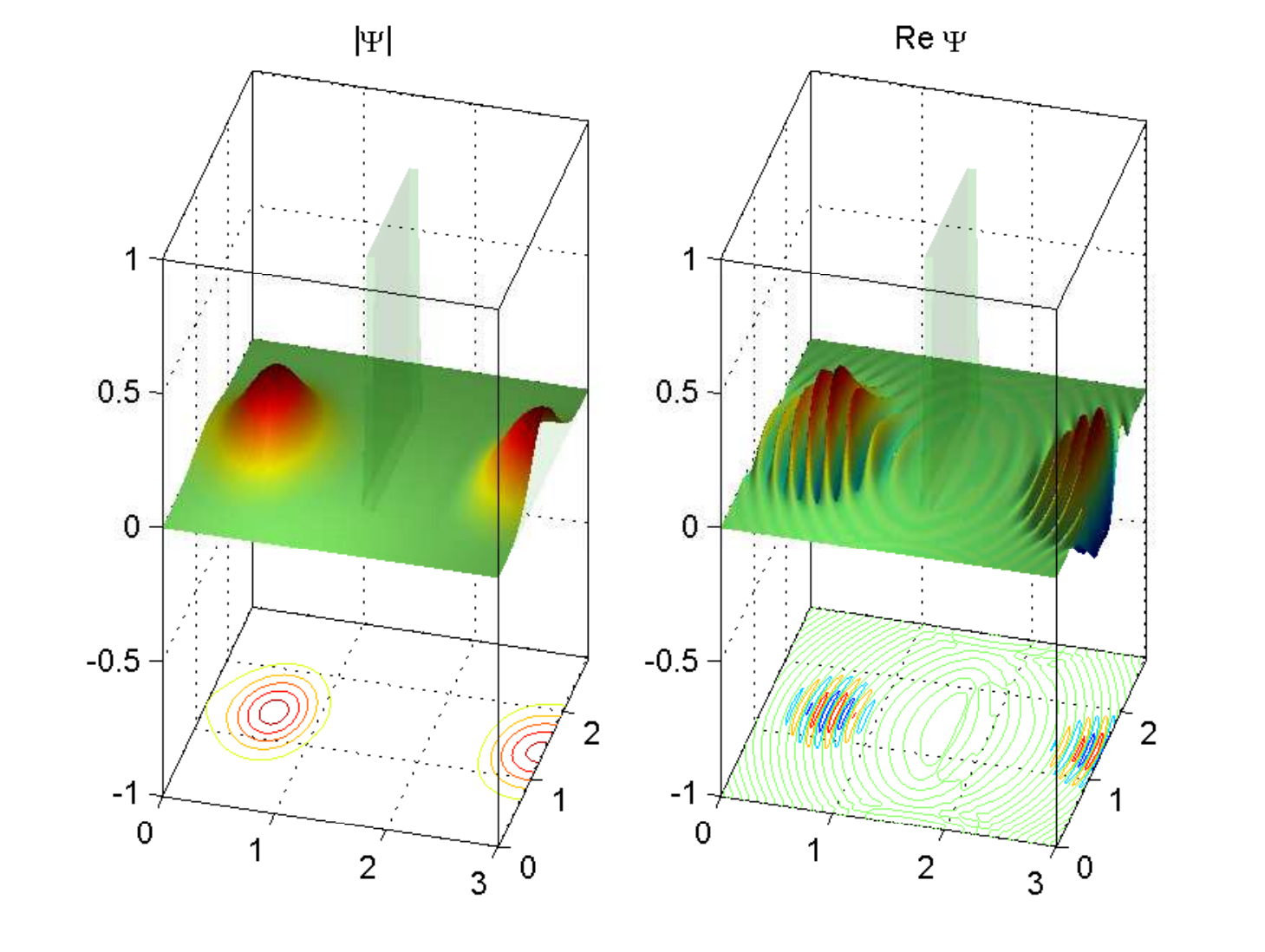}\vspace{0cm}\\%
\centerline{\small{$m=510$}}\\
\end{multicols}
\caption{\small{Example B. The modulus and the real part of the numerical solution $\Psi^m$, $m=204,264,360$ and $510$, for $(J,K,M)=(300,64,600)$}}
\label{SSP:EX22a:B:Solution}

\end{figure}

\smallskip\par\textbf{Example B.}
Next we consider example B from \cite{DZZ13} and take the rectangular potential (the barrier)
\[
 V(x,y)=
 \begin{cases}
Q &\text{for}\ \ (x,y)\in (a,b)\times (c,d)\\[1mm]
0 &\text{otherwise}
\end{cases},\ \ Q>0
\]
depending both on $x$ and $y$.
Specifically, we take $(a,b)\times (c,d)=(1.6,1.7)\times (0.7,2.1)$ and the barrier height $Q=1500$.
This barrier is discontinuous and thus more complicated than the previous one from the numerical point of view due to non-smoothness of the corresponding exact solution.
\par We choose $(X,Y)=(3,2.8)$, then $(a,b)\times (c,d)\subset \Omega_{X,Y}$ and $\psi_G$ is small enough outside $\bar{\Omega}_{X,Y}$.
Let also $T=t_M=0.027$.
\par Notice carefully that we take $J$ and $K$ so that the points $a$, $b$, $c$ and $d$ belong to the corresponding meshes in $x$ and $y$ and exploit the averaged mesh potential
\begin{equation}
 V_{hjk}=\begin{cases}
 V(x_j,y_k)   &\text{for}\ \ x_j\neq a,b\ \ \text{and}\ \ y_k\neq c,d\\[1mm]
 Q/2 &\text{for}\ \ x_j=a,b\ \ \text{but}\ \ y_k\neq c,d,
     \ \text{or for}\ \ y_k=c,d\ \ \text{but}\ \ x_j\neq a,b\\[1mm]
 Q/4 &\text{for}\ \ (x_j,y_k)=(a,c),(a,d),(b,c),(b,d)
\end{cases}
\label{num21}
\end{equation}
for any $j$ and $k$. We put $\tilde{V}=0$ and $\Delta V=V_h$.

\begin{figure}[ht]\centerline{
   \includegraphics[width=0.5\linewidth]{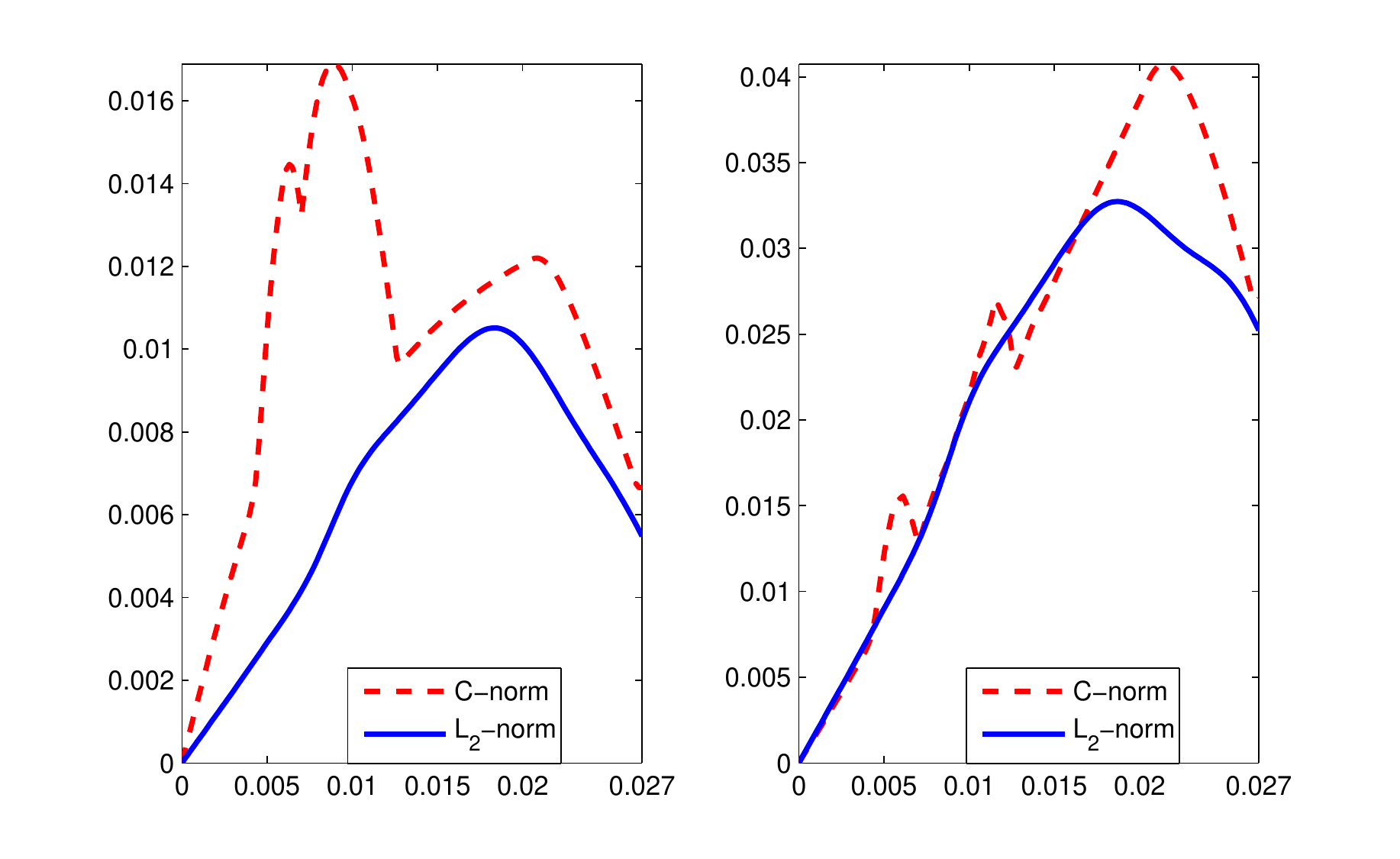}}
\caption{\small{Example B.
The absolute and relative differences in $C$ and $L^2$ norms between the numerical solutions for
$(J,K,M)=(300,64,600)$ and $(4800,256,2400)$, in dependence with time}
\label{SSP:EX22a:B:Errors}}
\end{figure}
\par The numerical solution $\Psi^{m}$ is computed for $(J,K,M)=(300,64,600)$
so that $h_x=10^{-2}$, $h_y=4.375\cdot10^{-2}$ and $\tau=4.5\cdot 10^{-5}$. Its modulus and real part are shown on Figure \ref{SSP:EX22a:B:Solution}, for the time moments $t_m=m\tau$, $m=204, 264, 360$ and $510$.
Once again the wave package is separated by the barrier (but now more abruptly) into two comparable reflected and transmitted parts moving in opposite $x$-directions and leaving the computational domain.
\par To justify the choice of the mesh with not so large $J$ (that is 4 time less than in \cite{DZZ13}) and $K$, on Figure \ref{SSP:EX22a:B:Errors} we give the absolute and relative differences in $C$ and $L^2$ space mesh norms between the numerical solutions for $(J,K,M)=(300,64,600)$ and $(4800,256,2400)$ that are 16 or 4 times larger, in dependence with time. Once again the level of differences is rather low.
\par To study the error behavior in more detail, we once again apply the approach described above.
In Table \ref{SSP:tab:EX22a:M}, the quantities
\begin{gather}
 \rho_{x,\bar{J},\bar{K},\bar{M},2^\ell},\ \ R_{x,\bar{J},\bar{K},\bar{M},2^\ell} \ \
 \text{for}\ \ (\bar{J},\bar{K},\bar{M})=(4800,256,2400)\ \ \text{and}\ \ \ell=5,4,3,2,1,
\label{num23}\\[1mm]
 \rho_{y,\bar{J},\bar{K},\bar{M},2^\ell},\ \ R_{y,\bar{J},\bar{K},\bar{M},2^\ell} \ \
 \text{for}\ \ (\bar{J},\bar{K},\bar{M})=(4800,512,2400)\ \ \text{and}\ \ \ell=5,4,3,2,1,
\nonumber \\[1mm]
 \rho_{t,\bar{J},\bar{K},\bar{M},2^\ell},\ \ R_{t,\bar{J},\bar{K},\bar{M},2^\ell} \ \
 \text{for}\ \ (\bar{J},\bar{K},\bar{M})=(4800,256,4800)\ \ \text{and}\ \ \ell=5,4,3,2,1
\nonumber
\end{gather}
(see \eqref{num3}-\eqref{num7}) are sequentially presented for the uniform in time and both $C$ and $L^2$ in space mesh norms.
We have also checked that quantities \eqref{num23} remain almost unchanged for less $\bar{M}=1200$.
\par Comparing the data with \eqref{num13}-\eqref{num15}, we conclude that now $\alpha\approx 2$ and $2<\beta\lesssim 3$ are less than in Example \textbf{A} whereas once again $\gamma\approx 2$. Note that the behavior of the approximations to $\beta$ is more chaotic than to $\alpha$ and $\gamma$ in both examples.
\begin{table}
\begin{center}
\begin{tabular}{|r|c|c|c|c|}
  \hline
   $J$ & $E_C$  & $E_{L^2}$  & $R_C$     & $R_{L^2}$  \\
  \hline
   150 & $0.663$$E$$-1$  &$0.286$$E$$-1$    & -       & -    \\
  300  & $0.118$$E$$-1$  &$0.390$$E$$-2$    & 5.62    & 7.33 \\
  600  & $0.270$$E$$-2$  &$0.907$$E$$-3$    & 4.38    & 4.30 \\
  1200 & $0.628$$E$$-3$  &$0.216$$E$$-3$    & 4.29    & 4.19 \\
  2400 & $0.125$$E$$-3$  &$0.435$$E$$-4$    & 5.03    & 4.99 \\
  \hline
\end{tabular}
\end{center}

\begin{center}
\begin{tabular}{|r|c|c|c|c|}
  \hline
    $K$ & $E_C$ &          $E_{L^2}$   &            $R_C$ & $R_{L^2}$  \\
  \hline
  16  & $0.720$$E$$-1$ & $0.348$$E$$-1$  &   -               &                -  \\
  32  & $0.440$$E$$-2$ & $0.320$$E$$-2$  & 16.36             &             10.87 \\
  64  & $0.611$$E$$-3$ & $0.379$$E$$-3$  & \hphantom{1}7.21  &  \hphantom{1}8.44 \\
  128 & $0.958$$E$$-4$ & $0.292$$E$$-4$  & \hphantom{1}6.37  &             12.90 \\
  256 & $0.223$$E$$-4$&  $0.378$$E$$-5$  & \hphantom{1}4.29  & \hphantom{1}7.79  \\
  \hline
\end{tabular}
\end{center}

\begin{center}
\begin{tabular}{|r|c|c|c|c|}
  \hline
   $M$ & $E_C$ &          $E_{L^2}$   &            $R_C$ & $R_{L^2}$  \\
  \hline
  150  & $0.172$\hphantom{$E$$-1$} & $0.917$$E$$-1$  &   -  &   -    \\
  300  & $0.430$$E$$-1$            & $0.233$$E$$-1$  & 3.97 &  3.93  \\
  600  & $0.106$$E$$-1$            & $0.580$$E$$-2$  & 4.37 &  4.02  \\
  1200 & $0.250$$E$$-2$            & $0.140$$E$$-2$  & 4.24 &  4.14  \\
  2400 & $0.506$$E$$-3$            & $0.274$$E$$-3$  & 4.94 &  5.11  \\
  \hline
\end{tabular}
\end{center}
\caption{Example B.
The change in numerical solution in maximum in time $C$ and $L^2$ space norms for redoubling
$J$ (and the reference $(\bar{J},\bar{K},\bar{M})=(4800,256,2400)$),
or $K$ (and the reference $(\bar{J},\bar{K},\bar{M})=(4800,512,2400)$),
or $M$ (and the reference $(\bar{J},\bar{K},\bar{M})=(4800,256,4800)$)}
\label{SSP:tab:EX22a:M}
\end{table}
\par
\begin{figure}[ht]
\begin{multicols}{2}
    \includegraphics[scale=0.37]{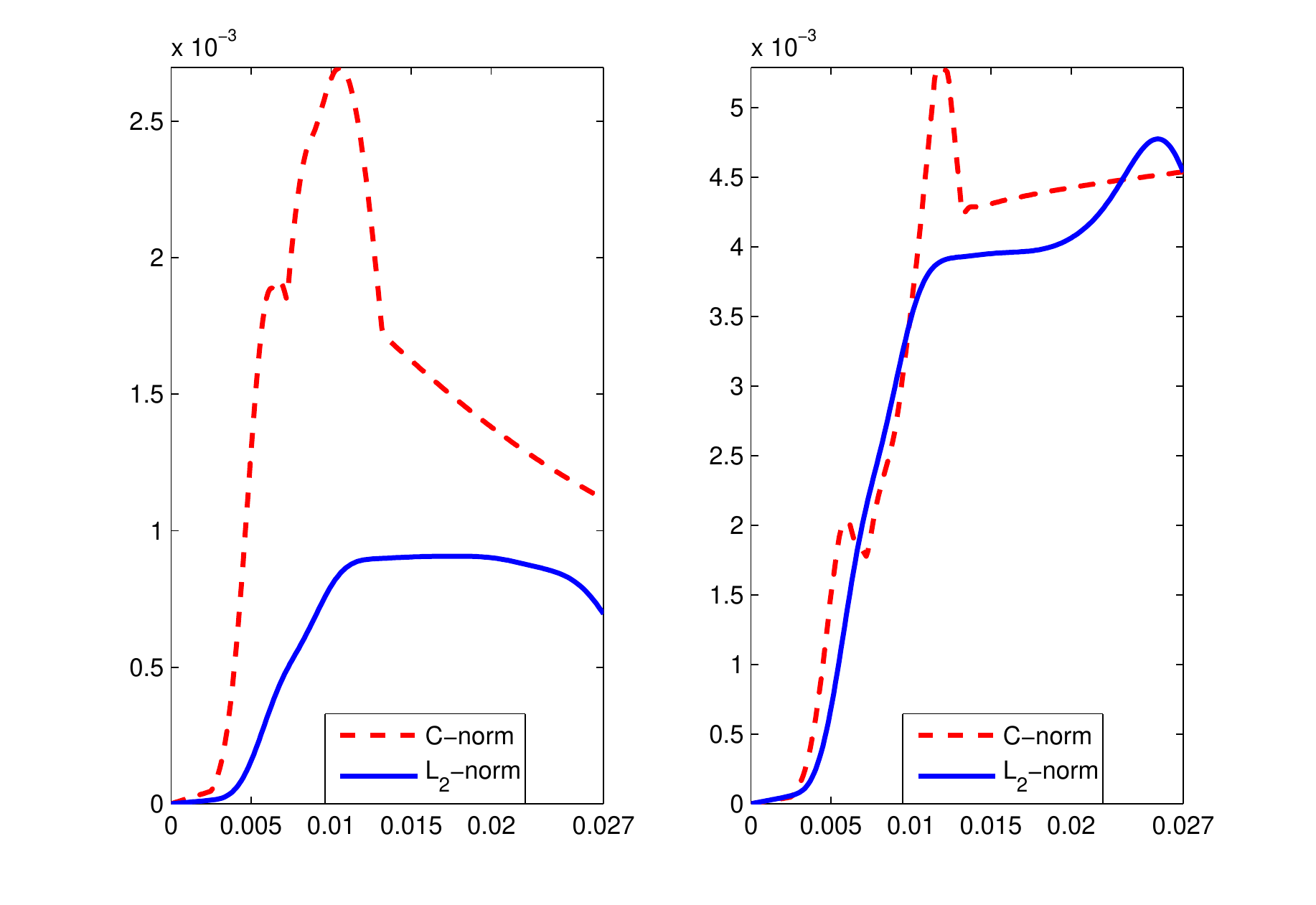}\vspace{0cm}\\
\caption{\small{Example B.
The absolute and relative differences in $C$ and $L^2$ norms between the numerical solutions for
$J=600$ and 4800, while
$(K,M)=(256,2400)$,
in dependence with time}}
\label{SSP:EX22a:B:Solution_J}
    \includegraphics[scale=0.4]{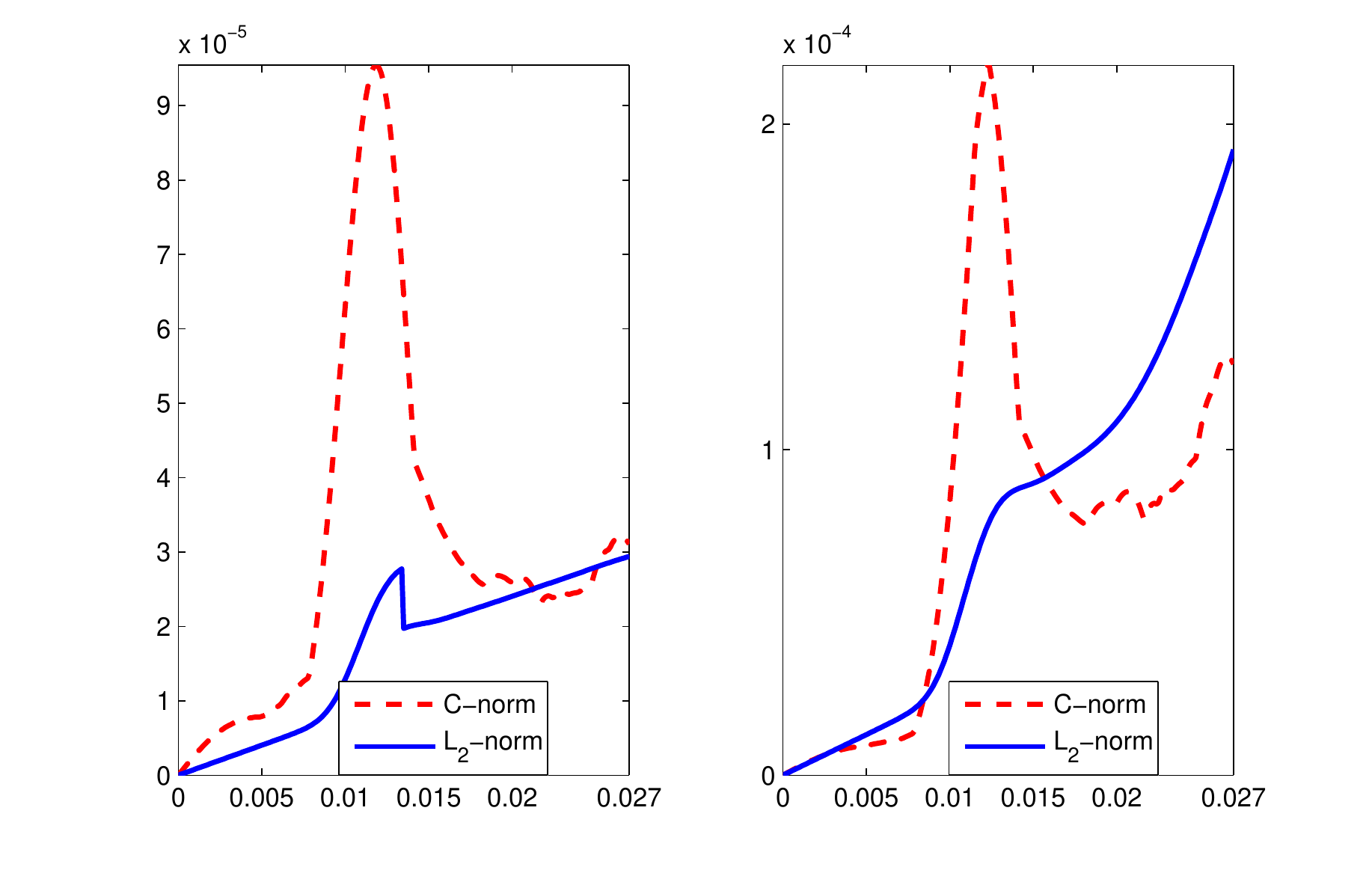}\vspace{0cm}\\
\caption{\small{Example B.
The absolute and relative differences in $C$ and $L^2$ norms between the numerical solutions for
$K=128$ and 512, while
$(J,M)=(4800,2400)$,
in dependence with time}}
\label{SSP:EX22a:B:Solution_K}
\end{multicols}
\end{figure}

 In addition, Figures \ref{SSP:EX22a:B:Solution_J} and \ref{SSP:EX22a:B:Solution_K} present the typical behavior of $\Psi_{\bar{J},\bar{K},\bar{M}}-\Psi_{\bar{J}/2^\ell,\bar{K},\bar{M}}$ and $\Psi_{\bar{J},\bar{K},\bar{M}}-\Psi_{\bar{J},\bar{K}/2^\ell,\bar{M}}$
both in $C$ and $L^2$ space norms, in absolute and relative forms, in dependence with time, for above $(\bar{J},\bar{K},\bar{M})$ and $\ell=3$ and $2$ respectively.

\par We also computed the differences between the numerical solutions using the averaged potential \eqref{num21} and the non-averaged one:
$E_C=0.724$$E$$-1$ and $E_{L^2}=0.330$$E$$-1$ for $(J,M,K)=(600,64,600)$ (compare with Figure \ref{SSP:EX22a:B:Errors})
and only twice (not four times) less
$E_C=0.367$$E$$-1$ and $E_{L^2}=0.166$$E$$-1$ for double $(J,M,K)=(1200,128,1200)$. Thus exploiting of the non-averaged potential reduces the orders of convergence significantly
(we omit more details).
\par Comparing the above error analysis with somewhat different one given in \cite{DZZ13}, we see the advantages of the present higher order in space splitting method allowing to exploit coarser space meshes even for discontinuous potentials and oscillating in space and time solutions that is not so obvious a priori.
\par We also notice that Figures \ref{SSP:EX22b:B:Solution2} and \ref{SSP:EX22a:B:Solution} for the numerical solutions together with Figures \ref{SSP:EX22b:B:Errors} and \ref{SSP:EX22a:B:Errors} for the corresponding relative errors (which decrease while the waves cross the artificial boundaries)
exhibit the absence of the spurious reflections from the artificial left and right boundaries; this is due to setting the discrete TBCs there.
\smallskip\par
\textbf{Acknowledgments}
\smallskip\par
This work is supported by
the National Research University Higher School of Economics' Academic Fund Program, project No. 13-09-0124 and by the Russian Foundation for Basic Research, project No. 12-01-90008-Bel.

\end{document}